%From saltman@mail.ma.utexas.edu  Mon Feb 26 15:33:12 2001
%% VERSION:  2/26/01
\magnification=\magstep1
\hsize=6.3truein
\input amssym.def
\input amssym.tex
%%%%%%%%%%%%  Macros necessary for this file  %%%%%%%%%%%%%5
\def\scong{{\scriptstyle\|}\lower.2ex\hbox{$\wr$}}
\def\Z{{\Bbb Z}}

\def\Ext{\mathop{\rm Ext}\nolimits}

\def\Tr{\mathop{\rm Tr}\nolimits}
\def\End{\mathop{\rm End}\nolimits}

\def\Res{\mathop{\rm Res}\nolimits}

\def\Hom{\mathop{\rm Hom}\nolimits}
\def\Ind{\mathop{\rm Ind}\nolimits}

\def\Gal{\mathop{\rm Gal}\nolimits}
\def\Br{\mathop{\rm Br}\nolimits}

\def\diagram{\def\normalbaselines{\baselineskip21pt\lineskip1pt
	\lineskiplimit0pt}}
\def\rtimes{\mathop{\times\!\!{\raise.2ex\hbox{$\scriptscriptstyle|$}}}
	\nolimits} 
\def\proof{\noindent{\it Proof.}\quad}
\def\blackbox{\hbox{\vrule width6pt height7pt depth1pt}} 
\outer\def\Demo #1. #2\par{\medbreak\noindent {\it#1.\enspace}
	{\rm#2}\par\ifdim\lastskip<\medskipamount\removelastskip
	\penalty55\medskip\fi}
\def\qed{~\hfill\blackbox\medskip}
\overfullrule=0pt
\def\Br{\mathop{\rm Br}\nolimits}

\def\P{{\Bbb P}}

\def\hangbox to #1 #2{\vskip1pt\hangindent #1\noindent \hbox to #1{#2}$\!\!$}

\pageno=0
\footline{\ifnum\pageno=0\hfill\else\hss\tenrm\folio\hss\fi}
\topinsert\vskip1.8truecm\endinsert
\centerline{INVARIANT FIELDS OF SYMPLECTIC AND ORTHOGONAL GROUPS}
\vskip6pt
$${\vbox{\halign{\hfil\hbox{#}\hfil\qquad&\hfil\hbox{#}\hfil\cr
$$David J. Saltman$^*$\cr
Department of Mathematics\cr
The University of Texas\cr
Austin, Texas 78712\cr}}}$$
\vskip16pt
\vskip3pt
{\narrower\smallskip\noindent
{\bf Abstract} The projective orthogonal and symplectic groups  
$PO_n(F)$ and 
\break 
$PSp_n(F)$ have a natural action on the 
$F$ vector space $V' = M_n(F) \oplus \ldots \oplus M_n(F)$. 
Here we assume $F$ is an infinite field of characteristic not 2. 
If we assume there is more than one summand in $V'$, 
then the invariant fields $F(V')^{PO_n}$ and $F(V')^{PSp_n}$ 
are natural objects. They are, for example, the centers 
of generic algebras with the appropriate kind of involution. 
This paper considers the rationality properties of these fields, 
in the case $1,2$ or $4$ are the highest powers of 2 that 
divide $n$. We derive rationality when $n$ is odd, or when 2 is the highest 
power, and stable rationality when 4 is the highest power. 
In a companion paper [ST]
joint with Tignol, we prove retract rationality when 
8 is the highest power of 2 dividing $n$. Back in this paper, 
along the way, we consider two generic ways of forcing 
a Brauer class to be in the image of restriction.
\bigskip

\medskip

\smallskip}

\footnote{}{*The author is grateful for support under NSF grant DMS-9970213. 
The author would also like to acknowledge the hospitality 
of the Schroedinger Institute during part of the preparation 
of this paper.}

\vfill\eject
\leftline{Introduction}
\medskip
In this paper $F$ will always be an infinite field of 
characteristic not 2. Set $V = M_n(F) \oplus \ldots \oplus M_n(F)$ 
to be a sum of $r \geq 2$ copies of $n \times n$ matrices. 
This $V$ is a representation of the projective linear group 
$PGL_n = GL(F)/F^*$ induced by conjugation acting diagonally. 
The invariant field $F(V)^{PGL_n}$ has become an algebraic 
and geometric object of remarkable importance, partly 
because $F(V)^{PGL_n}$ is the center $Z(F,n,r)$ 
of the so called generic division algebra $UD(F,n,r)$. 
Procesi observed that $Z(F,n,r)$ is rational 
(i.e. purely transcendental) over $Z(F,n,2)$. 
For this reason, it makes sense to only consider 
the case $V = M_n(F) \oplus M_n(F)$. 

In addition, there are important subgroups $\cal G$ 
of $PGL_n$ where the invariant field 
$F(V)^{\cal G}$ is also of significance. 
In particular, we will be interested in ${\cal G} = 
PO_n$ the projective orthogonal group, and when 
$n$ is even, ${\cal G} = PSp_n$ the projective 
symplectic group. Note that these groups will be defined 
in detail in section one. By [P2] p. 377-378 
and [R] p. 184, $F(V)^{\cal G}$ is the center of 
the generic central simple algebra with orthogonal 
respectively symplectic involution. 

There is a rich theory of central simple algebras 
with involution, which make it immediately clear that 
$F(V)^{\cal G}/F$ is, for both the ${\cal G}$ above, 
retract rational over $F$ when $m$ is odd and $n = m$, $2m$, or $4m$. 
In a companion paper [ST], joint with J.-P. Tignol, we will 
prove retract rationality in the case $n = 8m$.
In the case $n = 2^rm$ for $r > 3$, 
retract rationality is equivalent to questions about 
central simple algebras that, as of now, have no answer. 
But in the $r = m$, $2m$, $4m$ cases, one can ask 
whether $F(V)^{\cal G}$ has stronger properties, for example, 
that this field is rational or stably rational over $F$. 
Recall that $K/F$ is stably rational if and only 
if there is a field $L \supset K$ such that $L/F$ and 
$L/K$ are rational. 

In the case $n = 4m$, we show in section five 
that $F(V)^{\cal G}/F$ is stably rational, 
a not surprising result. It is somewhat 
more difficult to see that $F(V)^{\cal G}/F$ is 
rational in the cases $n = m$ and particularly when 
$n = 2m$. The main part of this paper is concerned 
with proving this last fact. Along the way 
we introduce 2 ``generic'' ways of forcing a 
Brauer class to be in the image of restriction. 

For the reader's convenience, let us summarize the known results in: 

\proclaim Theorem. Let $F$ be an infinite field of 
characteristic not 2. 
Let $V = M_n(F) \oplus M_n(F)$ with the standard action 
of $PGL_n$. View $PO_n, PSp_n \subset PGL_n$. Then 
\smallskip
a) If $n$ is odd, $F(V)^{PO_n}/F$ is rational. 
\smallskip
b) If $n = 2m$ and $m$ is odd, then $F(V)^{PSp_n}/F$ and 
$F(V)^{PO_n}/F$ are both rational. 
\smallskip
c) If $n = 4m$ and $m$ is odd, then $F(V)^{PSp_n}/F$ and 
$F(V)^{PO_n}/F$ are both stably rational. 
\smallskip
d) If $n = 8m$ and $m$ is odd, then $F(V)^{PSp_n}/F$ and 
$F(V)^{PO_n}/F$ are retract rational. 

Note that part a) is 1.2 below, part b) 
is 4.1 and 4.2 below, and part c) is 5.1 below. 
Part d) is proved in the above mentioned paper [ST], joint work 
with J.-P. Tignol. Let us, for completeness sake, 
record the general version of the above result 
(with a sketch of a proof). 
Recall that for an algebraic group ${\cal G}$, 
a {\bf good} representation (always algebraic 
representation) is one where a generic 
point has trivial stabilizer. 

\proclaim Corollary. Let ${\cal G} = PO_n$ for all $n$ or 
$PSp_n$ for $n$ even. Let $m$ be an odd integer. 
\smallskip 
a) Let $V' = M_n(F) \oplus \ldots \oplus M_n(F)$ 
be a direct sum of at least two copies of $M_n(F)$ 
with the usual action by ${\cal G}$. If 
$n = m$ or $2m$, then $F(V')^{\cal G}$ is rational 
over $F$. If $n = 4m$, $F(V')^{\cal G}$ is stably rational 
over $F$. If $n = 8m$, then $F(V')^{\cal G}$ is retract rational 
over $F$. 
\smallskip
b) Let $V''$ be a good representation of ${\cal G}$ 
and $F$ algebraically closed of characteristic 0. 
If $n = m$, $2m$, or $4m$, then $F(V'')^{\cal G}$ is stably 
rational over $F$. If $n = 8m$, then $F(V'')^{\cal G}$ 
is retract rational over $F$. 

\proof Of course, 
we have not defined retract rational 
here. We refer the reader to [S] or 
[LN] p. 77 for the definition 
and basic properties. To prove part a),   
note that a proof parallel to that of A.3 shows that 
$F(V')^{\cal G}/F(V)^{\cal G}$ is rational. 
(Alternatively, $F(V')^{\cal G}$ is the center of the generic central simple 
algebra in the appropriate number of variables and with the 
appropriate involution. This is known to be rational over the 
same object with 2 variables, i.e., $F(V)^{\cal G}$.) 
In part b), the assumption on $F$ is only necessary 
to apply Bogomolov's no name lemma, so $F(V'')^{\cal G}$ 
is stably isomorphic to $F(V)^{\cal G}$. This proves b).~\qed
 
To begin the paper proper, let us define some notation and 
recall some constructions. 
$F$ will always be our ground field. 
If ${\cal G}$ is an algebraic group, 
an $F$ representation of ${\cal G}$ is a finite dimensional $F$ 
vector space $V$ and a  
map ${\cal G} \to GL_F(V)$ which is a homomorphism as a map 
between algebraic groups. 
If $\phi: K \to L$ is an embedding of fields, 
and $V$ is a $K$ vector space, we write $V \otimes_{\phi} L$ 
to mean the tensor product $V \otimes_K L$ where $L$ is a $K$ 
vector space via $\phi$. If $A/K$, $B/K'$ are algebras over $K,K'$ 
respectively, and $K' \subset K$, we abbreviate 
$A \otimes_K (B \otimes_{K'} K)$ as $A \otimes_K B$.

Suppose $K,L$ are fields regular and separably generated over $F$. 
Then $K \otimes_F L$ is a domain. The field of fractions, $q(K \otimes_F L)$,  
of $K \otimes_F L$ we call the join of $K,L$ over $F$. 
Suppose next that $K/F$ is a finite separable extension, 
and $L/K$ is separably generated. Let $K' \supset K \supset F$ 
be the Galois closure of $K/F$, and let $H = \Gal(K'/K)$, 
$G = \Gal(K'/F)$. Choose coset representatives $\sigma_i$, 
of $H$ in $G$ where $i = 1, \ldots, r$. 
That is, assume $G = \cup_{i=1}^r \sigma_iH$ is a disjoint 
union. We can also view $\sigma_i$ as an embedding of $K$ in 
$K'$, and we write $\sigma_i(L) = L \otimes_{\sigma_i} K'$. 
Since we may assume $\sigma_1 = 1$, $\sigma_1(L) =  K' \otimes_K L$ 
has an action by $H$ via the action on $K'$. 
If we set $T$ to be the field join $L\sigma_2(L)\ldots\sigma_r(L)$ 
over $K'$, 
then $T$ has an natural $G$ action extending that on $K'$. 
The fixed field $T^G$ will be written $\Tr_{K/F}(L)$ 
and is called the transfer. Note that if $L$ is the field 
of rational functions of a $K$ variety $V$, then $\Tr_{K/F}(L)$ 
is the field of rational functions of the transfer variety of $V$ to $F$. 
Also note that if $L/K$ is rational, then $\Tr_{K/F}(L)/F$ 
is rational. 

Suppose $G$ is a finite group and $L$ is a field with a 
(perhaps trivial) $G$ action. Let $M$ be a lattice over $G$. 
That is, let $M$ be a finite generated $\Z[G]$ module which 
is free as an abelian group. The group ring is, 
as a ring, just the Laurent polynomial ring 
$L[x_1,x_1^{-1},\ldots,x_r,x_r^{-1}]$ and so $L[M]$ is a domain. The action 
of $G$ on $L$ and $M$ induces an action on $L[M]$, 
which in turn induces an action on the field of 
fractions $q(L[M]) = L(M)$. 

We will also require twistings of the above construction. 
Suppose $\alpha \in \Ext(M,L^*)$. That is, suppose 
$\alpha$ is an extension $0 \to L^* \to M' \to M \to 0$. 
Then there is an induced action of $G$ on $L[M]$ 
such that $L[M]^* \cong M'$. We write $L_{\alpha}[M]$ 
to mean $L[M]$ with this $\alpha$ twisted action. 
We write $L_{\alpha}(M)$ to mean $q(L_{\alpha}[M])$ 
with the extended action. In this paper 
we will make significant use of group 
cohomology, to which we refer the reader to 
[B] as a good general reference. 

Finally, we recall a bit about central simple 
algebras and the Severi-Brauer variety. 
Let $A/K$ be a central simple algebra by which we mean 
that $A$ is finite dimensional over its center $K$ 
and simple. The dimension of $A/K$ is always a square, 
$n^2$, and we call $n$ the degree of $A$. Of course 
such an $A$ always has the form of matrices $M_r(D)$ 
over a division algebra $D/K$. We say $A/K$, $B/K$ 
are Brauer equivalent if they are matrices over the same 
division algebra, or equivalently if $M_r(A) \cong M_s(B)$ 
for some $r,s$. The equivalence classes under this relationship 
form, of course, the Brauer group of $K$ which we write as $\Br(K)$. 
Recall that the product is induced by tensor product 
and the inverse of the class of $A/K$ is the class of 
the opposite algebra $A^{\circ}$. 
If $L \supset K$, then $A \to A \otimes_K L$ induces a 
group homomorphism $\Br(K) \to \Br(L)$ we call restriction.
We also use the word restriction to refer to the 
map $A \to A \otimes_K L$ on algebras (and not just classes). 

If $L/K$ is $G$ Galois, and $\gamma \in H^2(G,L^*)$, 
then we can form the crossed product $\Delta(L/K,G,\gamma)$.  
By this we mean the algebra $\oplus_{g \in G} \, Lu_{g}$ 
where $u_gx = \sigma(x)u_g$ for all $x \in L$, and 
$u_gu_h = c(g,h)u_{gh}$ for a 2 cocycle $c(g,h)$ in $\gamma$. 
The crossed product induces an isomorphism $H^2(G,L^*) 
\to \Br(L/K)$, where $\Br(L/K)$ is the kernel of the restriction 
map $\Br(K) \to \Br(L)$. 
Some of the properties of central simple algebras 
which we need are contained in the classical: 

\proclaim Theorem 0.1. (e.g. [J] p. 226f, [LN] p.34) 
a) Suppose $A/K$, $B/K$ are 
central simple 
algebras of degree $n$ and are also Brauer equivalent. 
Then $A \cong B$ over $K$. 
\smallskip
b) Every Brauer class contains a unique division algebra, which is the 
member of the class of minimal degree. 
\smallskip
c) Suppose $A/K$ 
has degree $n$, $r \geq 1$ and $(r,n)$ is the gcd. Then the Brauer class 
of $A \otimes_K \ldots \otimes_K A$ ($r$ times) contains 
an element of degree dividing $n/(r,n)$. 

In particular in any Brauer 
equivalence class we can talk about {\bf the} central simple 
algebra of degree $n$, assuming such exists. 

Let us recall and record the useful result of Endo-Miyata 
([EM] or e.g. [LN] p. 82). 

\proclaim Proposition 0.2. Let $G$ be a finite group. 
Let $L/K$ be a $G$ Galois extension and $V$ an $L$ vector space 
with a semilinear $G$ action. Then $L(V)^G/K$ is a rational 
extension. As a consequence, if $P$ is a permutation $G$ 
lattice, then $L(P)^G/K$ is rational. As another consequence, 
if $V,W$ are faithful $F$ representations of $G$,  
then $F(V \oplus W)^G/F(V)^G$ is rational. In particular, 
$F(V)^G$ and $F(W)^G$ are stably isomorphic. 

Suppose $A/K$ is central simple of degree $n$. 
The Severi-Brauer variety $SB(A)$ (e.g. [LN] p.89 
for what follows) 
can be defined as the variety of dimension $n$ 
right ideals of $A$, realized as a closed subvariety 
of the Grassmann variety of subspaces of $A$ of dimension $n$. 
If $A = \End_K(V)$, then $SB(A)$ is isomorphic to the 
projective space $\P_K(V)$. It follows that $SB(A)$ 
is irreducible and has a field of fractions 
we write as $K(A)$. Note that $K(A)$ is the generic splitting 
field of $A$ introduced by Amitsur. In particular, 
$A \otimes_K K(A) \cong M_n(K(A))$. 
If $A$ is split, $SB(A)$ is projective space so  
$K(A)$ is rational over $K$. 
As a matter of notation, if $A/K'$ is central simple 
and $K' \subset K$, we will abbreviate $K(A \otimes_{K'} K)$ 
as $K(A)$. 
We quote two results 
about the birational isomorphisms of Severi-Brauer 
varieties. 

\proclaim Theorem 0.3. Suppose $A/K$ is central simple 
of degree $n$. 
\smallskip
a) (Amitsur e.g. [LN] 13.29): $K(A) \cong K(A^{\circ})$. 
\smallskip
b) (Tregub [T]): Suppose $n$ is odd. Then 
$K(A) \cong K(B)$, where $B$ is the element of the Brauer 
class of $A \otimes_K A$ of degree $n$. 

As a tool in what follows, we need a different description 
of $K(A)$. To begin this discussion, 
let $G \supset H$ be a finite group and subgroup, and 
$K/F$ a separable $G-H$ extension. That is, if $L/F$ is the 
Galois closure of $K/F$, then $G$ is the Galois group of 
$L/F$ and $H$ is the subgroup corresponding to $L/K$. 
Form the $G$ lattice $\Z[G/H]$. Recall that this is the 
lattice with $\Z$ basis $\{u_{gH} | gH \in G/H\}$ 
such that $g'(u_{gH}) = u_{g'gH}$. In $\Z[G/H]$ form the 
sublattice $I \subset \Z[G/H]$ generated by all $u_{gH} - u_{g'H}$. 
Note that there is an exact sequence 
$0 \to I \to \Z[G/H] \to \Z \to 0$. 
The boundary map $\delta: \Z = H^0(G,\Z) \to H^1(G,I)$ 
defines a 1 cohomology class $\alpha \in H^1(G,I)$ 
we call the canonical class. One can compute that $\alpha$ 
has order $[G:H]$, generates $H^1(G,I)$, and the restriction 
$\alpha|_H = 0$. 

\proclaim Lemma 0.4. a) Let $M$ be any $G$ module. 
There is a map $d: \Ext(I,M) \to H^2(G,M)$, natural in $M$, 
whose image is precisely those element of $H^2(G,M)$ 
that split when restricted to $H$. 
This map is one to one if $H^1(H,M) = 0$. 
\smallskip
b) If $\beta: 0 \to M \to M' \to I \to 0$ corresponds 
to $\beta \in \Ext(I,M)$, then $d(\beta) = 
\delta(\alpha)$ where $\delta$ is the boundary in the long exact 
sequence associated to $\beta$. 

\proof  
There is an exact sequence 
$$\Ext^1(\Z[G/H], M) \to 
\Ext^1(I,M) \to \Ext^2(Z,M) \to \Ext^2(\Z[G/H], M).$$ 
Now $\Ext^i(\Z[G/H], M) = H^i(G, \Hom(\Z[G/H], M)) = $
(by Shapiro's Lemma e.g. [B] p.73) $ = H^i(H,M)$. 
Similarly, $\Ext^i(\Z, M) = H^i(G,\Hom(\Z,M)) = H^i(G,M)$. 
Thus the above exact sequence reads 
$$0 \to H^1(H,M) \to \Ext^1(I,M) \to H^2(G,M) \to H^2(H,M).$$
The map $H^2(G,M) \to H^2(H,M)$ is easily seen to be restriction, 
proving a). Part b) is an easy computation.~\qed 

The above lemma is sort of a ``splitting module'' 
construction. Assume $H^1(H,M) = 0$ and $\gamma \in H^2(G,M)$ 
is split when restricted to $H$. Then $\gamma$ defines 
a unique $\beta \in \Ext(I,M)$ with $d(\beta) = \gamma$. 
Let $0 \to M \to M' \to I \to 0$ correspond to $\beta$. 
Since $\beta$ obviously splits in $M'$, $\gamma$ 
splits in $M'$. Suppose $f: M \to N$ is any 
morphism of $G$ modules with $H^1(H,N) = 0$. Then clearly 
$f$ extends to $M'$ if and only if $f_*(\gamma) = 0$. 

The above discussion makes it natural to apply 0.4 to the 
case $M = L^*$, where $L$ is a field with faithful 
$G$ action. Of course, $H^1(H,L^*) = 0$ by Hilbert's theorem 90. 
If $K = L^H$ and $F = L^G$, then the set of $\gamma \in H^2(G,L^*)$ 
that restrict to 0 on $H$ is precisely the relative Brauer 
group $\Br(K/F)$. If $\gamma \in \Br(K/F)$, let $A/F$ 
be the corresponding central simple algebra with maximal subfield 
$K$ and $\beta: 0 \to L^* \to M' \to I \to 0$ the corresponding extension. 
Let $L_{\beta}(I)$ be as defined above.  

\proclaim Theorem 0.5. $L_{\beta}(I)^G$ is the Amitsur generic 
splitting field of $A$. That is, $L_{\beta}(I)$ 
is the function field of the Severi-Brauer variety of $A$. 

\Demo Remark. This is well known. When $H = 1$ a proof appears in 
[LN] p. 95 except that there is an exponent error 
changing $A$ to $A^{\circ}$. 
This error turns out to be unimportant because of Amitsur's 
result 0.1 a). 
We sketch a proof because a precise reference, in this generality, 
and in anything like this language, is hard to come by. 

\proof 
As is well known, $K \otimes_F A^{\circ} \cong \End_K(A)$ 
where $A$ is a $K$ vector space by the left action and $A^{\circ}$ 
acts on $A$ by right multiplication. Thus $L \otimes_F A^{\circ} 
\cong \End_L(L \otimes_K A)$. The $G$ action on $L \otimes_F A^{\circ}$ 
induces a $G$ action on $\P_L(L \otimes_K A)$ and the Severi-Brauer 
variety of $A^{\circ}$ is the variety given by descent. 

Viewing $A^{\circ} \subset \End_L(L \otimes_K A) \subset 
\End_F(L \otimes_K A)$, let $B$ be the centralizer of $A^{\circ}$
in $\End_F(L \otimes_K A)$. Since the centralizer of $L$ in $B$  
is the centralizer of $L \otimes_K A^{\circ}$ in $\End_F(L \otimes_K A)$, this 
centralizer is $L$. Thus $B$ is Brauer equivalent to $A$ with maximal 
subfield $L$. That is, $B = \Delta(L/F,G,\gamma)$. 
We can choose $u_g \in B$ such that $u_gx = g(x)u_g$ 
for all $x \in L$. Any such $u_g$ must form a basis 
of $B$ over $L$ and $u_gu_{g'} = c(g,g')u_{gg'}$ 
where $c(g,g')$ is a cocycle in $\gamma$ (see for example the proof of 7.2 
in [LN]). 
If $\beta \in B$ and $w \in L \otimes_K A$ we write 
$\beta \cdot w$ for the action of $\beta$ on $w$. 
If $g \in G$ and 
$L\alpha$ is a line in $L \otimes_K A$, i.e. a point in $\P_L(L \otimes_K A)$, 
then it is easy to check that $g(L\alpha) = u_g \cdot {L\alpha}$ 
and this, of course, is independent of the choice of $u_g$ 
(see [LN] p. 94). 

We can choose the $u_h$ for $h \in H$ 
such that $u_h \cdot (x \otimes a) = h(x) \otimes a$ for all 
$x \otimes a \in L \otimes_K A$. 
Then for any coset $gH \subset G$, we can write  
$v_{gH}$ as the element $u_g \cdot (1 \otimes 1)$. Since the $u_g$ 
span $B$, it is clear that the $v_{gH}$ are an 
$L$ basis of $L \otimes_K A$ and we use this basis to define 
projective coordinates for $\P_L(L \otimes_K A)$. 
Specifically, let the line $L(\sum_{gH} a_{gH}v_{gH})$ 
have projective coordinates $(a_{gH})_{gH \in G/H}$. 
Note that 
$$u_{g'} \cdot L(\sum_{gH} a_{gH}v_{gH}) = 
L(\sum_{gH} g'(a_{gH})c(g',g)v_{g'gH}).$$ 

Define the rational functions $d_{gH}$ for $gH \not= H$ 
by $d_{gH}(L(\sum_{gH} a_{gH}v_{gH})) = a_{gH}/a_H$, 
so the $d_{gH}$ are a transcendence basis for 
$L(\P_L(L \otimes_K A))$. Finally, 
$$(g'(d_{g''H}))(L(\sum_{gH} a_{gH}v_{gH})) 
= g'( d_{g''H}(u_{g'^{-1}} \cdot L(\sum_{gH} a_{gH}v_{gH})))= $$
$$g'(d_{g''H}(L(\sum_{gH} g'^{-1}(a_{gH})c(g'^{-1},g)v_{g'^{-1}gH}))) = $$
Using the definition of $d_{g''H}$ this is: 
$$(a_{g'g''H}/a_{g'H})g'(c(g'^{-1},g'g'')/c(g'^{-1},g')) = 
(a_{g'g''H}/a_{g'H})c(g',g'')^{-1}.$$ 
Thus $g'(d_{g''H}) = (d_{g'g''H}/d_{g'H})c(g',g'')^{-1}$. 
Replacing $A^{\circ}$ by $A$ and $c(g',g'')^{-1}$ by $c(g',g'')$ 
we have the result.~\qed 

It is convenient to use 0.5 to reprove the following 
result. 

\proclaim Corollary 0.6. Let 
$Z = Z(F,n)$ be the center of the generic division algebra 
$UD = UD(F,n)$. The generic splitting field $Z(UD)$ is rational 
over $F$. 

\proof 
We recall that $Z(F,n) = F(M)^{S_n}$ where there is an 
exact sequence $0 \to M \to M' \to I[S_n/S_{n-1}] \to 0$ 
and $M'$ is a permutation lattice with 
$M' = M'' \oplus \Z[S_n/S_{n-1}]$. 
The cocycle $\gamma \in H^2(S_n,M)$ defines the Brauer group 
element given by $UD$. By 0.5, $Z(UD) = F(M')^{S_n}$. 
But $F(\Z[S_n/S_{n-1}])^{S_n}$ is rational over $F$ 
since $F(\Z[S_n/S_{n-1}]) = F(x_1,\ldots,x_n)$ with the usual 
action. By Endo-Miyata (0.2), $F(M')^{S_n}/F(\Z[S_n/S_{n-1}])^{S_n}$ 
is rational.~\qed

We need to make an observation similar to 0.4, but involving 
cohomology of degree one higher. Let $H \subset G$ and $I$ 
be as in 0.4. For convenience, assume 
$G = H \cup HgH$. It immediately follows that $u_H - u_{gH}$ 
generates $I$. Choose $H' \subset H \cap gHg^{-1}$. 
Then there is surjection $\Z[G/H'] \to I$ defined 
by sending the canonical $H'$ generator to $u_H - u_{gH}$. 
Define $J$ to be the kernel of this map. Let $M$ be a $G$ 
module. Assume $\beta \in H^2(H,M)$. Denote by 
$g(\beta)$ the induced element of $H^2(gHg^{-1}, M)$, 
and $\beta - g(\beta)$ the obvious induced element 
of $H^2(H \cap gHg^{-1}, M)$. 

\proclaim Proposition 0.7. Let $M$  be a $G$ module  
such that $H^1(H',M) = 0$. 
\smallskip
a) Let $A \subset H^2(H,M)$ 
be the subgroup consisting of those $\beta$ such 
that the restriction 
of $\beta - g(\beta)$ to $H'$ is 0. 
For each $\beta \in A$ , there is an extension 
$0 \to M \to M' \to J \to 0$, corresponding to 
$\alpha \in \Ext(J,M)$ such that $\alpha = 0$ if and 
only if $\beta$ is in the image of $H^2(G,M)$. 
This correspondence is functorial in $M$. 
\smallskip
b) Suppose further that $H^3(G,M) \to H^3(H,M)$ 
is injective. Then this correspondence is onto 
$\Ext(J,M)$. 
\smallskip
c) Suppose $0 \to M \to M' \to J \to 0$ corresponds
to $\beta\! \in\! H^2(H,M)/\Res(H^2(G,M))$. 
Then the image of $\beta$ generates the kernel of 
$$H^2(H,M)/\Res(H^2(G,M) + \delta(H^1(H,J))) \to H^2(H,M')/\Res(H^2(G,M').$$ 

\proof We begin with a). 
As before, we have the exact sequence $\Ext(\Z[G/H'],M) \to \Ext(J,M) 
\to \Ext^2(I,M) \to \Ext^2(\Z[G/H'],M)$. In addition, 
we have the exact sequence 
$\Ext^2(\Z,M) \to \Ext^2(\Z[G/H],M) \to \Ext^2(I,M) \to \Ext^3(\Z,M) \to 
%\break 
\Ext^3(\Z[G/H],M)$. 
Since $\Ext(\Z[G/H'], M) = H^1(H',M)$, 
we have $\Ext(J,M) \to \Ext^2(I,M)$ is injective. 
In a similar way, 
$\Ext^2(\Z[G/H'], M) \cong H^2(H',M)$, 
%\break 
$\Ext^i(\Z[G/H], M) = H^i(H,M)$, and 
$\Ext^i(\Z,M) \cong H^i(G,M)$. 
Just as in 0.4, 
the induced map $H^i(G,M) \to H^i(H,M)$ is restriction. 
Thus we have exact sequences $H^2(G,M) \to H^2(H,M) \to 
\Ext^2(I,M) \to H^3(G,M) \to H^3(H,M)$ and 
$0 \to \Ext(J,M) \to \Ext^2(I,M) \to H^2(H',M)$. 
The elements of $H^2(H,M)$ which define elements 
of $\Ext(J,M)$ are precisely the elements in the kernel of 
$H^2(H,M) \to \Ext^2(I,M) \to H^2(H',M)$. If $A'$ 
is this kernel, then the injectivity of $\Ext(J,M) \to \Ext^2(I,M)$ 
shows that a) holds for $A'$. Thus to prove a) it suffices to observe:
 
\proclaim Lemma 0.8. The composition $h: H^2(H,M) \to \Ext^2(I,M) 
\to \Ext^2(\Z[G/H'],M) = H^2(H',M)$ is the map taking $\beta$ 
to the restriction $(\beta - g(\beta))|_{H'}$. 

\proof
Let $f: \Z[G/H'] \to \Z[G/H]$ be the induced map taking 
the canonical generator to $u_1 - u_{g(1)}$. Then $f$ 
induces $f': \Hom(\Z[G/H],M) \to \Hom(\Z[G/H'],M)$. $f'$, 
in turn, induces $f'^*: H^2(G,\Hom(\Z[G/H],M)) \to H^2(G,\Hom(\Z[G/H'],M))$. 
The map $h$ we are analyzing is the composition of 
$\phi_{H'}f'^*\phi_H^{-1}$ 
where $\phi_H: H^2(G,\Hom(\Z[G/H],M)) \to H^2(H,M)$ 
and $\phi_{H'}:H^2(G,\Hom(\Z[G/H'],M) \cong H^2(H',M)$ 
are the isomorphisms from Shapiro's 
Lemma ([B] p. 73 and p. 80). Note that as an $H$ module, 
$\Z[G/H] = \Z \oplus \Z[H/H \cap gHg^{-1}]$. 
In particular, $H^2(H,\Hom(\Z[G/H],M)) = 
H^2(H,M) \oplus H^2(H \cap gHg^{-1},M))$. 
Suppose $\phi_H^{-1}(\beta_H) = \beta_G \in H^2(G,\Hom(\Z(G/H),M))$. 
The restriction $\beta_G|_H$ is easily seen to 
be 
$$(\beta_H, \beta') \in H^2(H,M) \oplus H^2(H \cap gHg^{-1},M) = 
H^2(H,\Hom(\Z[G/H],M)).$$ 
Here $\beta'$ is the restriction 
of $g(\beta_H)$ to $H \cap gHg^{-1}$. 

The Shapiro Lemma ([B] p. 80) 
isomorphism $\phi_H: H^2(G,\Hom(\Z[G/H'],M)) \cong H^2(H',M)$ 
is a composition $H^2(G,\Hom(\Z[G/H'],M)) \to 
H^2(H',\Hom(\Z[G/H'],M)) \to H^2(H',M)$ where the 
first map is restriction and the second map, call it $\rho_{H'}^*$, 
is induced by the 
canonical $H'$ morphism $\rho_{H'}: \Hom(\Z[G/H'],M) \to M$. 
Thus in computing $h(\beta_H) = \phi_{H'}(f'^*(\beta_G)$, 
we can first restrict $\beta_G$ to $H'$. 
That is, $\phi_{H'}(f'^*(\beta_G) = {\rho'_{H'}}^*(f'^*|_H(\beta_G|_{H'}))$ 
where $f'^*|_{H'}$ is the induced map on $H'$ cohomology groups. 
Since $H' \subset H \cap gHg^{-1} \subset H$, the restriction, $\beta_G|_{H'}$, 
can be written $(\beta_H|_{H'}, g(\beta)|_{H'})$. 
View $f$, and hence $f'$ as $H'$ module maps (so $f'^*|_{H'}$ 
can be written simply as $f'^*$). It remains to compute 
$(\rho_{H'} \circ f'))^*(\beta_H|_{H'},g(\beta)|_{H'}$. 

Let $\rho: \Z \to \Z[G/H']$ 
be the $H'$ map sending $1 \in \Z$ to the canonical $H'$ fixed 
generator of $\Z[G/H']$. Then $\rho$ induces the 
$\rho_{H'}: \Hom(\Z[G/H'], M) \to 
\Hom(\Z,M) = M$ mentioned above. 
$\rho_{H'} \circ f': \Hom(\Z[G/H],M) \to \Hom(\Z,M) = M$ 
is induced by $f \circ\rho: \Z \to \Z[G/H] = \Z \oplus \Z[H/(H \cap gHg^{-1})]$ 
which we compute sends $1 \in \Z$ to $(1,-u_{gH}) \in 
\Z \oplus \Z[H/(H \cap gHg^{-1})]$. 
Now $u_{gH}$ is the canonical generator of $\Z[H/(H \cap gHg^{-1}]$, 
and this proves that $(\rho_{H'} \circ f'))^*(\beta_H|_{H'},g(\beta)|_{H'} 
= \beta_H|_{H'} - g(\beta)|_{H'}$. Thus a) is proven. 

We next consider b). 
By assumption, $H^3(G,M) \to H^3(H,M)$ has 0 kernel and so 
$H^2(H,M) \to \Ext^2(I, M)$ is surjective. If $\alpha \in \Ext(J,M)$, 
then $\alpha$ maps to, say,  $\alpha' \in \Ext^2(I,M)$ which is the 
image of some $\beta \in H^2(H,M)$. Since $\alpha'$ maps to 0 
in $H^2(H',M)$, $\beta \in A$ and b) is proven.

Finally we turn to showing part c) of 0.7. Let 
$\alpha \in \Ext(J,M)$ correspond to $0 \to M \to M' \to J \to 0$ 
and let $\beta \in H^2(H,M)$ be the preimage of $\alpha$. 
Of course, the canonical map $\Ext(J,M) \to \Ext(J,M')$ 
maps $\alpha$ to 0. It follows that the image of $\beta$ 
in $H^2(H,M')$ is also in the image of some $\delta \in H^2(G,M')$. 

We have the diagram: 
$$\diagram
\matrix{
H^2(G,M)&\longrightarrow&H^2(G,M')&\longrightarrow&H^2(G,J)&\longrightarrow&H^3(G,M)\cr
\downarrow&&\downarrow&&\downarrow&&\downarrow\cr
H^2(H,M)&\longrightarrow&H^2(H,M')&\longrightarrow&H^2(H,J)&\longrightarrow&H^3(H,M)\cr}$$
Suppose $\beta' \in H^2(H,M)$ is such that the image of 
$\beta'$ in $H^2(H,M')$ is also in the image of 
some $\delta' \in H^2(G,M')$. 
We need to show that $\beta'$ is a power of $\beta$ modulo 
$\Res(H^2(G,M)) + \delta(H^1(H,J))$. 

Let $\epsilon, \epsilon' \in H^2(G,J)$ be the image of 
$\delta, \delta'$ respectively. We claim it suffices to show 
that $\epsilon'$ is a power of $\epsilon$. To prove the claim, 
assume $\epsilon' = \epsilon^r$. Then $\delta' = \delta^r\chi$ 
where $\chi$ is the image of $\mu \in H^2(G,M)$. Modifying 
$\beta'$, $\delta'$ by the image of $\mu^{-1}$, we may assume 
$\delta' = \delta^r$. The result is now clear.

Thus we are reduced to showing $\epsilon'$ 
is a power of $\epsilon$. Tracing through the above diagram, 
$\epsilon'$ maps to 0 
in $H^2(H,J)$. Consider the exact sequence 
$H^1(G,I) \to H^2(G,J) \to H^2(G, \Z[G/H'])$. The map 
$H^2(H', \Z) \cong H^2(G,\Z[G/H']) \to H^2(H,\Z[G/H'])$ is injective 
because the inverse of the isomorphism is restriction to $H'$ 
followed by projection to $\Z$, and one can factor the restriction 
to $H'$ step  through restriction to $H$. 
It follows that $\epsilon'$ maps to 0 in $H^2(G, \Z[G/H'])$ 
and hence is in the image of $H^1(G,I)$. 
But tracing through the equivalences, it is easy but tedious 
to see that $\epsilon$ is the image of a generator 
of $H^1(G,I)$, proving c).~\qed 

Let us take the construction in 0.7 and apply it to fields. 
In 0.7, let $M$ have the form $L^*$ where $L/F$ is a $G$ 
Galois extension. Recall that an extension 
$\gamma: 0 \to L^* \to M' \to J \to 0$ defines a ``twisted'' 
action of $G$ on $L[J]$ and the field of fractions $L(J)$. 
If this extension is defined by $\beta \in H^2(H,L^*) = 
\Br(L/L^H)$, we write $L[J]$ and $L(J)$ with this twisted action 
as $L_{\beta}[J]$ and $L_{\beta}(J)$ respectively. 
Note that we can form $L_{\beta}(J)$ for any $\beta$ 
with $(\beta - g(\beta))|_{H'} = 0$. 
Also note that if $\beta$ is the image of $\alpha \in \Br(L^G)$ 
then since $L$ splits $\beta$, $L$ must split $\alpha$. 
That is, $\alpha$ defines an element in $H^2(G,L^*)$. 
Translating 0.7. we have:

\proclaim Proposition 0.10. 
The extension of $\beta$ to $\Br(L_{\beta}(J)^H)$ is in the image of 
$\Br(L_{\beta}(J)^G)$. If $\beta$ is in the image of 
$\Br(L^G)$, then $L_{\beta}[J] \cong L[J]$ 
the isomorphism preserving $G$ actions and $L[J]$ 
having untwisted action. 

\bigskip

\leftline{Section One: Reducing the finite group}
\medskip

We begin  by recalling the definition 
of the projective orthogonal and symplectic groups. 
Of course, $PGL_n(F) = GL_n(F)/ F^*$ 
is the quotient of $GL_n(F)$ modulo its center. 
The orthogonal group $O_n = O_n(F)$ 
is the subgroup of $GL_n(F)$ where $AA^T = I$ 
and $T$ refers to the transpose. 
We define $PO_n = PO_n(F)$ to be the image 
of $O_n(F)$ in $PGL_n(F)$. Note that this means 
$PO_n(F)$ is not necessarily the $F$ points of the 
corresponding algebraic group scheme. To avoid this technicality 
we would have to introduce $GO_n$, the group of so called 
similitudes (e.g. [K-T] p. 153). However, 
the subgroup $PO_n(F)$, as we have defined it, 
is Zariski dense in the $\bar F$ points (i.e. over the algebraic closure) 
and in considering invariant rings or fields this issue 
is therefore irrelevant. 
%\vfill\eject

Next we recall the definition of the symplectic group. 
The symplectic involution $J_1$ on $2 \times 2$ matrices 
is defined as 
$$\pmatrix{a&b\cr c&d\cr} \longrightarrow \pmatrix{d&-b\cr -c&a\cr}\ .$$
We identify $M_{2m}(F)$  with $M_2(F) \otimes_F M_m(F)$ 
so that the matrix idempotent $e_{11} \otimes e_{ii}$ is 
$e_{2i-1,2i-1}$ and $e_{22} \otimes e_{ii}$ is $e_{2i,2i}$. 
The symplectic involution $J_m$ on $M_{2m}$ can then be described as 
$J_1 \otimes T$ where $T$ is the transpose on $M_m(F)$. From now on, 
we write $J_m$ as $J$. 
Of course, $Sp_n$ is the group of matrices $A \in M_{2m}(F)$ 
such that $A^JA = AA^J = I_{2m}$. 
We define $PSp_n(F)$ to be the image of $Sp_n$ in $PGL_n(F)$, 
so $PSp_n(F) = Sp_n(F)/\{I,-I\}$.
Once again $PSp_n(F)$ are not the $F$ points of the group scheme, 
and once again it does not matter. 

Let $V = M_n(F) \oplus M_n(F)$ with the natural 
diagonal action of $PGL_n(F)$. Then 
$F(V)^{PGL_n}$ is the center $Z = Z(F,n,2)$ 
of the generic division algebra of degree $n$ in 2 
variables which we write as $UD = UD(F,n,2)$. 
As remarked before, combining [P2] p. 377-78 and [R] p. 184 
we have that 
$F(V)^{PO_n}$ and $F(V)^{PSp_n}$  are the centers of the 
generic algebras with appropriate involution. 
Combining this with [BS] 
p. 112, , we have the following. 

\proclaim Theorem 1.1. a) If $n$ is even, $F(V)^{PSp_n}/F(V)^{PGL_n} = 
F(V)^{PSp_n}/Z$ 
is the function field of the 
Severi Brauer variety of an algebra $A/Z$ 
where $A$ has degree $n(n-1)/2$ and is Brauer equivalent 
to $UD \otimes_Z UD$.  
\smallskip
b) $F(V)^{PO_n}/F(V)^{PGL_n} = 
F(V)^{PO_n}/Z$ 
is the function field of the Severi Brauer variety of an algebra $A/Z$ 
where $A$ has degree $n(n+1)/2$ and is Brauer equivalent 
to $UD \otimes_Z UD$.  

Recall that if $B/K$ is any central simple algebra, we write 
$K(B)$ to mean the function field of the Severi-Brauer variety 
of $B$. That is, $K(B)$ is the Amitsur generic splitting field 
of $B$. By e.g. [LN] 13.12, $K(M_r(B))$ is rational over $K(B)$. 
In particular, let $D/Z$ be the division algebra Brauer equivalent 
to $UD \otimes_Z  UD$, which it is not hard to see has degree $m = 
n/2$ or $n$ depending on whether $n$ is even or odd. 
One result we desire is now easy. 

\proclaim Theorem 1.2. Suppose $n$ is odd. 
Then $F(V)^{PO_n}$ is rational 
over $F$. 

\proof By the result 0.3 b) of Tregub, $Z(D) \cong Z(UD)$. 
By 0.6, $Z(UD)/F$ is rational.~\qed

We can now assume $n = 2m$ is even. 
Using 1.1, both $F(V)^{PO_n}$ and $F(V)^{PSp_n}$ 
are rational over $Z(D)$ of degree $n(n+1)/2 - m$ and $n(n-1)/2 - m$ 
respectively. In particular, to prove rationality or stable rationality 
results for $F(V)^{PSp_n}$ and $F(V)^{PO_n}$, it suffices 
to consider the case $PSp_n(F)$. 

Procesi showed that $Z(F,n ,2) = F(V)^{PGL_n}$ can be written 
as a multiplicative invariant field $F(M)^{S_n}$ 
where $M$ is a lattice whose definition we will recall 
later. Here $S_n$ is the symmetric group in $n$ letters, and so 
is the Weyl group of $PGL_n$ ([P] or e.g. [LN] p. 109). 
In this section we will recall 
the parallel argument for $PSp_n$, and so get 
$F(V)^{PSp_n}$ as a multiplicative invariant field 
of the Weyl group, $W$, of $PSp_n$. In addition, we will 
prove an analogue of 1.1 where $Z$ is replaced by $F(M)^W$. 
This is a significant improvement, as $W$ is smaller than 
$S_n$ and has an abelian normal subgroup we will use. 

To this end, let $T_{PGL} = T_{GL}/F^*$  where $T_{GL} 
\subset GL_n(F)$ is the group of diagonal matrices. 
Of course $T_{PGL}$ is a maximal torus of $PGL_n(F)$. 
Let $S_n \subset PGL_n$ 
be the group of permutation matrices, so $N_{PGL_n}  = 
T_{PGL_n}S_n$ is the normalizer of $T_{PGL_n}$ and 
the inclusion $S_n \subset N_{PGL_n}$ induces $S_n \cong 
N_{PGL_n}/T_{PGL_n}$. 

Let $T_{Sp} \subset Sp_n$ 
be the subgroup of $Sp_n$ of diagonal matrices. 
The condition defining $Sp_n$ implies that any such 
diagonal matrix has the form  
$(a,a^{-1},b,b^{-1},\ldots,c,c^{-1})$ down the diagonal. 
Set $T_{PSp} \subset PSp_n$ to be the image. 
$T_{PSp}$ is a maximal torus of $PSp_n$. 

Let $N_{PSp}$ 
be the normalizer of $T_{PSp}$ in $PSp_n$. Then $W = N_{PSp}/T_{PSp}$ 
is the Weyl group and can be described as follows. 
Writing $M_{2m}(F) =  M_2(F) \otimes M_m(F)$, consider 
$\tau_i = \tau \otimes e_{ii} + \sum_{j \not= i}I_2 \otimes e_{jj}$ 
where: 
$$\tau = \pmatrix{0&-1\cr 1&0\cr}$$
%%insert here 
Note that $\tau_i^2 \in T_{Sp}$, and that $\tau_i$ normalizes 
$T_{Sp}$. Let $A \subset W$  be generated by the $\tau_i$, 
so that $A \cong \Z/2\Z \oplus \ldots \oplus \Z/2\Z$ ($m$ times). 
Let $S_m \subset M_m(F)$ be the group of permutation matrices 
embedded in $M_{2m}(F)$ via $\sigma \to 1 \otimes \sigma$, 
and identifying $S_m$ with this image, $S_m$ normalizes 
$T_{Sp}$. One can show that $W \cong A \rtimes S_m$ where 
$S_m$ has the obvious permutation action on the $\tau_i$. 
Said differently, $W \subset S_{2m}$ is the subgroup 
preserving the partition ${\cal P} = \{\{1,2\},\{3,4\},\ldots,\{2m-1,2m\}\}$. 

Having the definition of $W$ in front of us, 
let us define some objects involving $W$. 
Let $p_i': A \to \Z/2\Z$ be defined by 
$p_i'(\tau_j) = 0$ for $j \not= i$ and $p_i'(\tau_i) = 
1 + 2\Z$. Let $S_{m-1} \subset S_m$ 
be the subgroup fixing both elements in $\{1,2\}$, and 
set $H = A \rtimes S_{m-1}$. 
Note that $H$ is the stabilizer in $W$ of the set $\{1,2\}$. 
By setting $p_1(S_{m-1}) = 0 + 2\Z$, $p_1'$ extends 
to $p_1: H \to \Z/2\Z$.

The embedding $PSp_n \subset PGL_{2m}(F)$ induces 
the commutative diagram:
$$\diagram
\matrix{
T_{Sp}&\subset&T_{GL}\cr
\downarrow&&\downarrow\cr
T_{PSp}&\subset&T_{PGL}\cr}\leqno{(1)}$$
We also have $N_{PSp} \subset N_{PGL}$. Since $T_{PGL} \cap N_{PSp} = T_{PSp}$, 
we have an induced embedding $W \subset S_{2m}$ which 
is precisely the one described above. 
Finally, note that we can define an intermediate group 
$N_{PSp} \subset N' \subset N_{PGL}$ such that $N' \supset T_{PGL}$ 
and $N'/T_{PGL} = W$. 

As a module over $PSp_n$, $M_n(F)$ is a direct sum of the spaces of 
symmetric and skew symmetric matrices. Let $M^-$ be the 
later. Of course, $M^-$ is the Lie algebra of $PSO_n$. 
Let $\Delta \subset M^-$ be the sub vector space 
of diagonal matrices, necessarily of the form 
$(a,-a,b,-b,\ldots,c,-c)$ down the diagonal. 
Note that $\Delta$ is preserved by $N_{PSp}$.

Since $M^{-}$ is the Lie algebra of $G$, 
there is a invariant Zariski open subset $U \subset M^-$ 
such that, for all $x \in U$, $Gx \cap \Delta \not= \emptyset$ 
and if $gx,g'x \in D$, then there is an $h \in N_{PSp}$ 
such that $hgx = g'x$. 

Let $V^-$ be the $PSp_n$ representation $M^- \oplus M_n(F)$. 
Then by A.1, 
we can write $F(V^{-})^{PSP_n} = F(\Delta \oplus M_n(F))^{N_{PSp}}$. 
Note that $T_{PGL}$ acts trivially on $\Delta$, 
so $V^- \subset V$ is preserved by $N'$ and the 
action of $N_{PSp}$ on $V^-$ extends to $N'$. 

Set $K = F(\Delta)$ and view $F(\Delta \oplus M_n(F))$ as $K(M_n(F))$. 
If the $y_{ij} \in K(M_n(F))$ correspond to the standard basis 
of $M_n(F)$, then $K(M_n(F))$ is the field of fractions 
of $K[Y']$ where $Y' \subset K(M_n(F))$ is the multiplicative 
subgroup generated by the $y_{ij}$. The elements of $Y'$ 
(i.e. monomials in the $y_{ij}$) in  
$K[Y']$ form a basis of $T_{PGL}$ and hence $T_{PSp}$ 
eigenvectors. We next evaluate $K[Y']^{T_{PGL}}$ and 
$K[Y']^{T_{PSp}}$. 

The character group $\Hom_F(T_{PGL},F^*)$ 
is a sublattice of $\Hom(T_{GL},F^*) =$\break $  \Z[S_{2m}/S_{2m-1}]$.  
More precisely, let $d_i: T_{GL} \to F^*$ 
be the projection on the $i^{th}$ diagonal entry. 
Then the $d_i$ form a basis of $\Hom(T_{GL},F^*)$ permuted 
by $S_{2m}$ in the obvious way. $\Hom(T_{PGL},F^*)$ 
is the sublattice generated by all elements of the 
form $d_i - d_j$. Note that $d_i - d_j$ is precisely 
the character associated to the eigenvector $y_{ij}$. 
Set $I_{2m} = \Hom(T_{PGL},F^*)$. Then we have a  surjective $S_{2m}$ 
morphism $Y' \to I$ associating each monomial to its character, 
and we set $Y$ to be the kernel. It is clear that $K[Y']^{T_{PGL}}$ 
is spanned by the $T_{PGL}$ fixed monomials and so 
$K[Y']^{T_{PGL}} = K[Y]$.  

To compute $K[Y']^{T_{PSp}}$ we begin by computing $M_{Sp} = \Hom(T_{Sp},F^*)$ 
as an image of $\Z[S_{2m}/S_{2m-1}]$ under restriction. 
{From} the structure 
of $T_{Sp}$ we have that the character $f_i = d_{2i} + d_{2i-1}$ 
is trivial on $T_{Sp}$. 
%In fact, 
It is clear that $M_{Sp} = 
\Z[S_{2m}/S_{2m-1}]/\Z[S_m/S_{m-1}]$  where $\Z[S_m/S_{m-1}]$ is the 
sublattice generated by the $f_i$. If we set $M_{PSp} \subset M_{Sp}$ 
to be $\Hom(T_{PSp}, F^*)$, it is clear from (1) that $M_{PSp} = 
I_{2m}/(\Z[S_m/S_{m-1}] \cap I_{2m})$. Furthermore, one can compute that 
$\Z[S_m/S_{m-1}] \cap I_{2m}$ is the lattice generated by all $f_i - f_j$ 
which we write as $I_m$. Note that $W$ acts naturally on all these 
lattices via its action on all the tori. Set $H' = W \cap S_{2m-1}$ 
so as a $W$ lattice, $\Z[S_{2m}/S_{2m-1}] = \Z[W/H']$. 
$A$ acts trivially on $\Z[S_m/S_{m-1}]$ and the 
induced action by $S_m = W/A$ is just the permutation action indicated. 
That is, $\Z[S_m/S_{m-1}]$ is the lattice $\Z[W/H]$ where 
$H \subset W$ is generated by $S_{m-1}$ and $A$. 
All together we have the following diagram of $W$ lattices 
corresponding to the diagram of tori (1): 
$$\diagram
\matrix{
\Z[W/H']/\Z[W/H]&\longleftarrow&\Z[W/H']\cr
\uparrow&&\uparrow\cr
I_{2m}/I_m&\longleftarrow&I_{2m}\cr}$$

Once again, there is a $W$ morphism 
$Y' \to I_{2m}/I_m$ taking each monomial to its 
$T_{PSp}$ character. Clearly this morphism is just the composition 
$Y'  \to I_{2m} \to I_{2m}/I_m$. Furthermore, if $Y_2$ 
is the kernel of $Y' \to I_{2m}/I_m$, then $K[Y']^{T_{PSp}} = 
K[Y_2]$. We have begun the proof of: 

\proclaim Proposition 1.3. 
a) $K[Y']^{T_{PGL}} = K[Y]$; $K[Y']^{T_{PSp}} = K[Y_2]$ 
\smallskip
b) $Y \subset Y_2$ and $Y_2/Y \cong I_m$ 
\smallskip
c) $K(Y')^{T_{PGL}} = K(Y)$; $K(Y')^{T_{PSp}} = K(Y_2)$

\proof 
We have already shown a), and b) follows from the fact that 
$Y_2$ is the inverse image of $I_m$ under the map $Y' \to I_{2m}$. 
As for c), this follows from the standard: 

\proclaim Lemma 1.4. 
Suppose ${\cal G} \subset GL_F(V')$ is an algebraic group 
with ${\cal G}_0 \subset {\cal G}$ the connected component 
of the identity. Let $U' \subset V'$ be a ${\cal G}$  invariant 
basic Zariski open subset (including the case $U' = V'$).  
If $\tau: {\cal G}_0 \to F^*$ is any character, 
assume there is a $u \in F[U']$ such that $\eta(u) = 
\tau(\eta)^{-1}u$ for all $\eta \in {\cal G}_0$. 
Then the field of fractions $q(F[U']^{\cal G})$ is 
$F(V')^{\cal G}$. 

\proof Write $F[U'] = F[V'](1/d)$. 
Since $F[V']$ is a UFD, 
$F[U']^*/F^*$ is a free abelian group with basis 
corresponding to the primes dividing $d$. 
For all $\eta \in {\cal G}_0$, $\eta(d) \in dF^*$ 
so ${\cal G}_0$ must permute the primes dividing 
$d$. That is, ${\cal G}_0$ acts trivially on $F[U']^*/F^*$.

Since $F[U']$ is a UFD, if $\alpha \in F(V')$ is $\cal G$ 
invariant, we may assume $\alpha = f/g$ where $f$, $g$ 
have no common factors. Since ${\cal G}/{\cal G}_0$ is finite, 
it suffices to find $f,g$ which are ${\cal G}_0$ invariant. 
If $\eta \in {\cal G}_0$, $g\eta(f) = \eta(g)f$. 
Since $f,g$ have no common factors, $\eta(f) = f\tau(\eta)$ 
and $\eta(g) = g\tau(g)$ where $\tau(g) \in F[U']^*$. 
Since ${\cal G}_0$ acts trivially on $F[U']^*/F^*$, 
the map ${\cal G}_0 \to F[U']^*/F^*$ induced by $\tau$ 
must be a homomorphism. Since ${\cal G}_0$ is connected, 
it follows that $\tau(g) \in F^*$ and $\tau: {\cal G}_0 \to F^*$ 
is a character. Choose $u \in F[U']$ as in the given. 
Then $\alpha = uf/ug$ and $uf,ug \in F[U']^{{\cal G}_0}$.~\qed

Having described the $T_{PSp}$ invariant field, the first part of the 
next proposition is clear. As for the rest, note that 
it says $F(V)^{PSp_n}$ is ``too big'' and the important 
information resides in a smaller field. Not only 
do we make it smaller by substituting $F(V^-)$ for 
$F(V)$, but we observe the ``$\Delta$'' part is also
irrelevant. 

\proclaim Proposition 1.5. $F(V^-)^{PSp_n} = K(Y_2)^W$. 
$F(V)^{PSp_n}/F(V^-)^{PSp_n}$ is rational of transcendence degree 
$n(n-1)/2$. $K(Y_2)^W/F(Y_2)^W$ is rational of degree $m$, 
so all together $F(V)^{PSp_n}/F(Y_2)^W$ is rational of degree 
$2m^2$. 

\proof  
$F(V)^{PSp_n} = F(M^+ \oplus V^-)^{PSp_n}$ where $M^+$ 
is the submodule of $M_n$ consisting of $J$ symmetric matrices. 
The second statement now follows from A.3. 
$K = F(\Delta)$ so $K(Y_2) = F(Y_2)(\Delta)$. Since $W$ acts 
linearly on $\Delta$, the rest follows from 
the result of Endo-Miyata (0.2). ~\qed

Of course, we have exact sequences $0 \to I_m \to \Z[W/H] \to \Z \to 0$ 
and $\beta_2: 0 \to Y \to Y_2 \to I_m \to 0$. By 0.4, this 
second sequence is associated to an element of $\gamma_2 \in H^2(W,Y)$ 
split by $H$. More precisely, it is associated with the 
element $\delta_2(\alpha_m)$ where $\alpha_m \in H^1(W,I_m)$ 
is the canonical generator and $\delta_2$ is the boundary 
of the long exact sequence associated to $\beta_2$. 
Of course $Y$ is also part of the $S_{2m}$ sequence 
$\beta: 0 \to Y \to Y' \to I_{2m} \to 0$ 
and this defines $\gamma \in H^2(S_{2m},Y)$ as 
$\delta(\alpha_{2m})$ where $\alpha_{2m} \in H^1(S_{2m},I_{2m})$ 
is the canonical generator and $\delta$ is the boundary associated with 
$\beta$. There is a $W$ embedding $I_m \to I_{2m}$ defined above 
and direct computation shows that the image of $\alpha_m$ 
is twice the restriction of $\alpha_{2m}$ to $W$. 
Thus $\gamma_2$ is twice the restriction of $\gamma$ to $W$. 

The standard argument, which parallels the one proving 1.5, 
shows that $Z(F,n,2) = F(V)^{PGL_n} = F(X \oplus Y)^{S_n}$ 
where $X = \Z[S_n/S_{n-1}]$. Another direct computation 
shows that $\gamma$ is the cocycle associated to the 
generic division algebra $UD(F,n,2)$.  Since $\Delta$ has rank $m$ 
and $X$ has rank $n = 2m$, $F(X \oplus Y)^W/F(Y)^W$ is rational 
of transcendence degree $2m$. It follows that $F(Y)^W$ can be thought 
of as the center of a generic division algebra with maximal subfield 
having $W$ as Galois group. This division algebra is described by the 
restriction of $\gamma$ to $W$ and we call it $D_{\gamma}$. 
Let $D_2$ be the division algebra of degree $m$ Brauer equivalent 
to $D_{\gamma} \otimes_{F(Y)^W} D_{\gamma}$. 

The sequence $0 \to Y \to Y_2 \to I_m \to 0$ induces 
the sequence $\beta': 0 \to F(Y)^* \to F(Y)^*Y_2 \to I_m \to 0$ 
and the field $F(Y_2)$ with its $W$ action can be  thought of as 
$F(Y)_{\beta'}(I_m)$. Thus by 0.5 we have the following. 

\proclaim Lemma 1.6. $F(Y_2)^W = F(Y)^W(D_2)$. 
That is, $F(Y_2)^W/F(Y)^W$ is the function field of the Severi 
Brauer variety of $D_2$. 

Set $U$ to be the 
$F$ representation $F[W/H']$ of $W$, 
where we recall that $H' = W \cap S_{2m-1}$.
With its $W$ action, $F(X)$ is just $F(U)$.  
In particular, by Endo-Miyata (0.2)) again, $F(X \oplus Y)^W/F(Y)^W$ 
is purely transcendental of degree $2m$. Finally, as we mentioned 
before, for any central simple algebra $A/L$, $L(M_r(A))$ 
is purely transcendental over $L(A)$. Using these facts, and adding the 
appropriate number of indeterminants to the fields involved, 
we have the following. 

\proclaim Theorem 1.7. Let $A$ be the central simple algebra 
of degree $n(n-1)/2$ in the Brauer class of 
$\Delta(F(Y)/F(Y)^W, W, \gamma^2)$. Let $A'$ 
be the central simple algebra in the same Brauer class 
of degree $m$. 
\smallskip
a) $F(V)^{PSp_n} \cong F(X \oplus Y)^W(A)$. 
\smallskip
b) $F(V)^{PSp_n}$ is rational over $F(Y)^W(A')$ of degree 
$2m^2$. 

Note that a) above is almost the same as 1.1 a) except the group 
$S_n = S_{2m}$ has been replaced by $W$. This is what is meant 
by ``reducing the group'' in the title of this section. 

\bigskip
\leftline{Section Two: More about lattices}
\medskip
We want to look further at the $W$ lattice  
$Y$ defined in section one. To begin with, we note that as a 
lattice over 
$W$, $Y'$ is the direct sum of two sublattices, we call 
$Y_D'$ and $Y_O'$. To be exact, $Y_D'$ is spanned by the 
$y_{ij}$ where $i,j$ are in the same element of the partition 
${\cal P} = \{\{1,2\},\ldots,\{2m-1,2m\}\}$, 
and $Y_O'$ is the span of the rest of the $y_{ij}$. 
Graphically, the $y_{ij}$ in $Y_D'$ correspond to matrix entries 
on the $2 \times 2$ matrix diagonal of the matrix $(y_{ij})$ 
and the $y_{ij} \in Y_O'$ 
correspond to the off diagonal entries. 
Also note that $Y'_O = \Z[W/H'']$ where $H''$ is the stabilizer 
of $y_{13}$. That is, $H''$ is $A' \rtimes S_{m-2}$ where 
$S_{m-2} \subset S_m$ is the stabilizer of $\{1,2,3,4\}$ and 
$A' \subset A$ is the kernel of the projection  
$p_1' \oplus p_2': A \to \Z/2\Z \oplus \Z/2\Z$. 

Consider the image of $Y_D'$ under the morphism $Y' \to I_{2m-1}$. 
Clearly this image is generated by $d_i - d_j$ where $i,j$ are in the 
same element of the partition ${\cal P}$. 
We wish to describe this image further. 
To this end, let $H = A \rtimes S_{m-1}$ be the subgroup 
fixing $d_1 + d_2$ just as in section one. 
Then the image in question is clearly isomorphic 
to $\Z^{-}[S_m/S_{m-1}] = \Ind_H^W(\Z^{-})$ where 
$\Z^{-} = \Z(d_1 - d_2)$. 
Of course $\Z^{-}$ is the rank one lattice associated to 
the homomorphism $p_1: H \to \{1,-1\} \cong \Z/2\Z$.   
Set $H' \subset H$ to be the kernel of this map, so 
$H' = W \cap S_{2m-1}$ as in section one. 

It is obvious that the cokernel 
of $\Z^{-}[S_m/S_{m-1}] \to I_{2m-1}$ is $I_m$ where 
$I_m$ fits into the sequence $0 \to I_m \to \Z[W/H] \to \Z$. 
All together we have the diagram:
%\vfill\eject
$$\diagram
\matrix{
&&0&&0&&0\cr
&&\uparrow&&\uparrow&&\uparrow\cr
0&\to&Y_O&\to&Y'_O&\to&I_m&\to&0\cr
&&\uparrow&&\uparrow&&\uparrow\cr
0&\to&Y&\to&Y'&\to&I_{2m-1}&\to&0\cr
&&\uparrow&&\uparrow&&\uparrow\cr
0&\to&Y_D&\to&Y_D'&\to&\Z^{-}[S_m/S_{m-1}]&\to&0\cr
&&\uparrow&&\uparrow&&\uparrow\cr
&&0&&0&&0\cr}$$

Let us consider further some of the lattices appearing above. 
Define $M \subset Y_D$ to be the $H$ sublattice 
generated by $y_{11}, y_{22}, y_{12} + y_{21}$. 
Of course, as an $H$ module, $M = \Z[H/H'] \oplus \Z$. 
It is then easy to see that $Y_D  = \Ind_H^W(M)$. 
In particular $Y_D$ is a permutation lattice. 
Thus $H^1(W',Y_D) = 0$ for any subgroup of $W' \subset W$.  
Shapiro's Lemma ([B] p. 73, 80) says 
$H^3(W,Y_D) \cong H^3(H,M)$ where the map is restriction 
to $H$ followed by projection onto $M$. If follows that 
$H^3(W,Y_D) \to H^3(H, Y_D)$ is injective. 

Note that $W = H \cup HgH$ where $g(1) = 3$ and $g(2) = 4$. 
Then $H \cap gHg^{-1}$ contains the group $H''$ defined above. 
We have  
$H^1(H'', Y_D)  = 0$. Since $Y_O' = \Z[W/H'']$ and the map 
$Y_O' \to I_m$ takes the $H''$ fixed generator to $g(1) - 1$, 
$Y_O$ 
is a lattice of the form of $J$ in 0.7. In all, 0.7 applies 
to the extension $0 \to Y_D \to Y \to Y_O \to 0$. 
We label this extension as $\beta'$, so 
if we set $L' = F(Y_D)$, then there is a $W$ 
action preserving field isomorphism $F(Y) \cong L'_{\beta'}(Y_O)$. 

Of course the above observation is most useful if we compute 
the element $\beta'' \in H^2(H,Y_D)$ associated to $\beta'$. 
Using 0.7 c) we can begin to do this by computing the kernel of 
$H^2(H,Y_D)/(\Res(H^2(W,Y_D)) + \delta(H^1(H,Y_O))) \to 
H^2(H,Y)/\Res(H^2(W,Y))$. 
Recall that $\Z^-[S_m/S_{m-1}]$ is spanned by all $d_i - d_j$ 
where $i,j$ are in the same element of ${\cal P}$. 
Let $\alpha \in H^1(H,\Z^-[S_m/S_{m-1}])$ be given by 
the cocycle $h \to h(d_1) - d_1$. Then clearly the image of 
$\alpha$ in $H^1(H,I_{2m-1})$ is the image of of the canonical 
generator of $H^2(W,I_{2m-1})$. Let $\beta \in H^2(H,Y_D)$ 
be the image of $\alpha$ under the boundary of the long exact sequence. 
Then the naturality of this boundary shows that 
the image of $\beta$ in $H^2(H,Y)$ is the restriction of the 
canonical element $\gamma \in H^2(W,Y)$. 

It will 
be convenient to describe this $\beta$ precisely. 
To begin with, $Y_D = \Ind_H^W(M)$ so $H^2(H,Y_D) = 
H^2(H,M) \oplus H^2(H \cap gHg^{-1}, g(M))$ (e.g. [B] p. 69). 
One can compute that $\beta = (\beta_H,0)$ where 
$\beta_H \in H^2(H,M)$  is the inflation, via $p_1$, of the 
nontrivial element 
of $H^2(\Z/2\Z, \Z(y_{12} + y_{21}))$. 

We will need to define a similar element in 
$H^2(H \cap gHg^{-1}, gM)$. Note that $g(M)$ can be 
viewed as the submodule 
spanned by $y_{33}$, $y_{44}$, and $y_{34} + y_{43}$. 
We define $p_2: gHg^{-1} \to \Z/2\Z$ to be the $g$ 
translation of $p_1$, and we note $p_2$ restricts to 
$p_2'$ on $A$. We also use $p_2$ to denote the restriction 
to $H \cap gHg^{-1}$. Then we set $\beta_{gH} \in 
H^2(H \cap gHg^{-1}, g(M)) 
\subset H^2(H,Y_D)$ to be the inflation via $p_2$ 
of the nontrivial element of $H^2(\Z/2\Z, \Z(y_{34} + y_{43}))$.

We will show $\beta''\! =\! \beta$ in some cases. 
To do this using 0.7 we compute 
$\delta(H^1(H,Y_O))\! \subset$ $ H^2(H,Y_D)$ which first requires 
computing $H^1(H,Y_O)$. 
But this last cohomology group 
is the image of $I_m^H$ under the boundary for the top horizontal 
sequence above. In fact, if $f_i \in \Z[S_m/S_{m-1}]$ 
is the image of $d_{2i-1}$ and $d_{2i}$ in $\Z[S_{2m}/S_{2m-1}]$, 
then $\eta' = (m-1)f_1 - \sum_{i \not= 1} f_i$ generates $I_m^H$.
Let $\eta$ be the image of $\eta'$  in $H^1(H,Y_O)$. 
Since we ultimately only need the case $m$ is odd we assume 
this and compute that $\eta$ has image $(0,\beta_{gH}) \in H^2(H,Y_D)$. 

We avoid further tedious cohomology details and just assert 
that $\beta$ is not in the image of $H^2(W,Y_D)$ 
but $(\beta_H,\beta_{gH})$ is. 
In particular $0 \to Y_D \to Y \to Y_O \to 0$ 
is not split. By 0.7 c), $\beta$ and $\beta''$ differ 
by an element of $H^2(W,Y_D) + \delta(H^1(H,Y_O))$. 
Since $\beta + \mu$ is in $\Res(H^2(W,Y_D))$, 
it follows that $\beta - \beta'' \in \Res(H^2(W,Y_D))$ 
and we might as well take $\beta'' = \beta$. 

\proclaim Theorem 2.1. 
The sequence $0 \to Y_D \to Y \to Y_O \to 0$ is associated, 
as in 0.7, with the element $\beta \in H^2(H,Y_D)$ defined 
above. 

In algebra terms, $\beta$ corresponds to the quaternion 
algebra $B \! =\!  \Delta(F(Y_D)^{H'}\! /F(Y_D)^H$, $y_{12}y_{21})$ 
over $F(Y_D)^H$. To see this recall that $H' = 
W \cap S_{2m-1}$ can also be described as the kernel 
of $p_1$. The gist of 2.1 is that the extension 
$F(Y)^W/F(Y_D)^W$ is a generic extension forcing $\beta$ 
to be in the image of $\Br(F(Y)^W)$. 
\bigskip
\vfill\eject 

\leftline{Section three: A cheaper way}
\medskip
In the last section, we described a generic way of forcing 
a Brauer group element to be in the image of restriction. 
In this section we describe a ``cheaper'' way, by which 
we mean a way that requires smaller transcendence degree. 
We will study its 
properties and study its connection with a 
construction involving generic matrices. 
But we approach this whole subject from the generic 
division algebra side, and in fact will start by considering 
a natural invariant theory problem involving generic matrices. 

To begin, we recall with slightly different emphasis an argument 
from [S1] (or [LN] p. 113). Let $D/F$ be a central simple algebra of degree 
$n$, and $UD = UD(F,r,s)$ the generic division algebra 
of degree $r$ in $s$ variables. We will also write 
$UD$ as $UD(F,P_1,\ldots,P_s)$ where the $P_i$ 
are the generic matrices generating $UD$. Let $Z = Z(F,r,s)$ 
be the center of $UD$ which we also write as $F(P_1,\ldots,P_s)$. 
Note that we avoid using $X_i$, $Y_i$ 
because in a future argument these $P_i$ will be $2 \times 2$ 
generic matrices and not $n \times n$. 

Recall that for $A/K$ any central simple algebra, $K(A)$ 
is the generic splitting field of $A$, or equivalently, 
$K(A)$ is the field of fractions of the Severi-Brauer variety 
defined by $A$. With this notation, set 
$K = Z(UD \otimes_K (D)^{\circ})$ where we ask 
the reader to recall our tensor product convention from the 
introduction. Let the $u_i$ be a basis of $D/F$.  
In [S1] or [LN] p. 113 we showed that: 

\proclaim Proposition 3.1. 
There is a canonical isomorphism $\phi: UD \otimes_Z K \to D \otimes_F K$. 
$K/F$ is rational with transcendence base $a_{ik}$ 
where $\phi(P_k) = \sum a_{ik}u_i$. 

By 3.1 we can think of $UD \otimes_Z K$ as $D(P_1, \ldots, P_s)$ 
where the $P_k = \sum a_{ik}u_i$ and $D(P_1,\ldots,P_s)$ has 
center $F(a_{ik})$. 
In particular this applies to the case where $D$ itself is a 
generic division algebra $UD(F',P,Q)$. 
We can use an old observation of Procesi 
to note: 

\proclaim Lemma 3.2. In the situation of 3.1, suppose 
$D = UD(F',P,Q)))$ so $F = Z(F',r,2)$. Then $K$ is the center of 
$UD(F',P,Q,P_1,\ldots,P_s) \cong UD(F,P,Q) \otimes_Z K$. 

\proof As above, let the $u_i$ be a basis 
of $UD(F',P,Q)$ over $F$. 
Procesi observed that $F'(P,Q,P_1,\ldots,P_s)$ is rational over $F'(P,Q)$ 
with transcendence basis the $y_{ik}$ where $P_k = \sum y_{ik}u_i$. 
That is, 3.2 is a direct consequence of 3.1.~\qed

Next we consider an invariant theory question. Let $S_m$ 
be the symmetric group, and let $UD(F,P_1,\ldots,P_m)$ 
be the generic division 
of degree $r$ in $m$ variables with center $F(P_1, \ldots,P_m)$. 
Abbreviate these algebras or fields as $UD(F,\vec P)$ and 
$F(\vec P)$ respectively. 
Then $S_m$ acts in the natural way on $UD(F, \vec P)$ and hence 
on $F(\vec P)$ by permuting the $P_i$'s. To be precise, we view 
$S_m$ as the bijections of $\{1,\ldots,m\}$ and 
$\sigma(P_i) = P_{\sigma(i)}$ for all $\sigma \in S_m$. 
Of course, the obvious question is the rationality 
of the invariant field $F(\vec P)^{S_m}$. It is not clear how to settle 
this, but it can be shown that this field is stably rational 
over the center of the appropriate generic division algebra. 
The next lemma is all we need in this direction, 
but see the remark for a bit more. 

\proclaim Lemma 3.3. Form the rational extension field 
$F(P,Q,\vec P) = F(P,Q,P_1,\ldots,P_m)$ where $S_m$ acts trivially on $P,Q$. 
Then 
$F(P,Q,\vec P)^{S_m}$ is rational over $F(P,Q)$. 
Similarly, the invariant 
field of 
$F(P,Q,P_1,\ldots,P_m,Q_1,\ldots,Q_m) = F(P,Q,\vec P,\vec Q)$ with the 
obvious $S_m$ action is rational over $F(P,Q)$. 

\Demo Remark. $UD(F,P_1,\ldots,P_m)^{S_m}$ is a division algebra 
of degree $r$ and so arguments like those of 3.3 show that 
that $F(P_1,\ldots,P_m,P,Q)^{S_m}$ is rational over 
$F(P_1,\ldots,P_m)^{S_m}$. It follows that $F(P_1,\ldots,P_m)^{S_m}$ 
is stably rational. 

\proof Write $F(P,Q,\vec P) = F(P,Q, y_{ik})$ where 
$P_k = \sum y_{ik}u_i$ and $u_i$ is a basis of $UD(F,P,Q)/F(P,Q)$. 
Clearly the action of $S_m$ on the $y_{ik}$ is just 
$\sigma(y_{ik}) = y_{i\sigma(k)}$. 
$F(P,Q,\vec P)^{S_m}$ is rational over $F(P,Q)$ 
via the usual fact about invariant fields of $S_m$ with permutation actions. 
The second sentence follows in the same way.~\qed

Set $L = F(P,Q,\vec P,\vec Q)$ with the given $S_m$ action, 
and now set $F_m = L^{S_m}$. We will use 3.2 to analyze the fields 
$L/F_m$ further. Set $L_k = F(P,Q,P_k,Q_k)$. It is clear that 
$L$ is the field join $L_1L_2\ldots{L_m}$ over $F(P,Q)$ 
with the obvious induced $S_m$ action. 
Set $UD = UD(F,P,Q)$ so $UD$ has center $F(P,Q)$. 
Set $UD_k = UD(F(P,Q),P_k,Q_k)$ and let $Z_k$ be the center 
of $UD_k$. By 3.2, 
$L_k = Z_k(UD_k \otimes_{Z_k} (UD)^{\circ}))$. 
Let $Z$ be the field join $Z_1\ldots{Z_m}$  with the obvious 
$S_m$ action. View $Z = Z_1\ldots{Z_m} \subset L_1\ldots{L_m} = L$ 
in the indicated way so that the $S_m$ actions are compatible. 

Denote by $S_{m-1} \subset S_m$ the stabilizer of $1$, and set $Z' =
Z^{S_{m-1}}$. We also set $F_m' = Z^{S_m}$. 
We can choose a set of left coset representatives $\sigma_k$ 
of $S_{m-1}$ in $S_m$ so that $\sigma_1$ is the identity 
and $\sigma_k(1) = k$. 
Then $Z_1 \subset Z'$ and we can set $UD_1' = UD_1 \otimes_{Z_1} Z'$, 
$UD' = UD \otimes_{F(X,Y)} Z'$, and $L' = Z'(UD_1' \otimes_{Z'} UD'^{\circ})$. 

\proclaim Theorem 3.4. $F_m$ is isomorphic to the transfer over $Z'/F'_m$ 
of the field extension $L'/Z'$. 

\proof The transfer is the field of fractions of the invariant ring 
of 
$$S = \bigotimes_{k=1}^m \sigma_k(Z \otimes_{Z'} L')$$ 
where $\sigma_k(Z \otimes_{Z'} L')$ refers to the $\sigma_k$ 
twist and the iterated tensor product is over $Z$. 
To prove the theorem, it suffices to find an $S_m$ 
invariant embedding $S \to L$ such that $L$ is the field 
of fractions of the image of $S$. This map is actually 
obvious, being just $\phi(l_1 \otimes \ldots \otimes l_m) =  
\prod_k \sigma_k(l_k)$. That the field of fractions of $\phi(S)$ 
is $L$ is clear, and that $\phi$ is injective can be seen by checking 
transcendence degrees.~\qed

Next we analyze $Z/F(P,Q)$ a bit more. We begin with $Z_1$ 
which is the center of the generic division algebra $UD(F(P,Q),n,2)$. 
As such, $Z_1$ has the form $F(P,Q)(M_r)^{S_r}$ where 
$M_r$ is the $S_r$ lattice described in section 
one. Form the 
wreath product group $W_r = (S_r \oplus \ldots \oplus S_r) \rtimes S_m$ 
where there are $m$ terms in the direct sum and the action 
of $S_m$ is the obvious one. Let $A \subset W_r$ be the $m$ fold 
direct sum $S_r \oplus \ldots \oplus S_r$ and 
$p_1': A \to S_r$ the projection 
on the first term. If we set $H = AS_{m-1} \subset W_r$, then $p_1'$ extends to 
$p_1: H \to S_r$ by setting $p_1(S_{m-1}) = 1$. Using $p_1$ 
we can view $M_r$ as an $H$ module and set $N = \Ind_{H}^G(M_r)$. 
We claim: 

\proclaim Theorem 3.5. $Z \cong F(P,Q)(N)^A$ with an isomorphism 
preserving the $S_m$ actions. Thus $F_m' \cong F(P,Q)(N)^{W_r}$. 

\proof 
Write $Z_1' = F(P,Q)(M_r)$ so $Z_1'^{S_r} = Z_1$. Thus, 
as fields (i.e. ignoring group actions), $Z'' = Z_1'\ldots{Z_m'} = 
F(P,Q)(M_r \oplus \ldots \oplus M_r)$ and this later field has an obvious 
action  by $A$ so that 
$Z = Z_1\ldots{Z_m} = F(P,Q)(M_r \oplus \ldots \oplus M_r)^A$. 
Checking the $S_m$ action makes it clear that $Z''$ 
has an induced $W_r$ action and that with respect to this action 
$Z'' = F(P,Q)(N)$. The theorem is now clear.~\qed 

If we set $r = 2$ the above picture will begin to look 
very familiar. $W_2$ is the Weyl group $W$, and $M_2$ 
is the lattice $M \oplus {\Z}x_1 \oplus {\Z}x_2$, where 
$M$ is as in the previous section. $N = 
\Ind_H^W(M \oplus {\Z}x_1 \oplus {\Z}x_2) 
= Y_D \oplus X$ where $X = {\Z}x_1 \oplus \ldots \oplus {\Z}_{2m}$.

We form the 
generic $2 \times 2$ generic division division algebra 
over $F(Y_D)^W$ and write its center, following the 
above conventions, as $F(Y_D)^W(P,Q)$. We claim:  
 
\proclaim Lemma 3.6. 
$F(Y_D)^W(P,Q)(y_1,\ldots,y_{2m}) = F(P,Q)(N)^W$  is the field 
$F_m'$ in the case $r = 2$. 

\proof 
This follows immediately because $F(P,Q)(N) = F(P,Q)(Y_D \oplus X)$ 
and $X$ is a permutation lattice.~\qed

Next we look at the quaternion algebra 
$\Delta(F(Y_D)^{H'}/F(Y_D)^H, y_{12}y_{21})$ from the previous section.  
Tracing through the definitions, we find that: 

\proclaim Lemma 3.7. 
$\Delta(F(Y_D)^{H'}/F(Y_D)^H, y_{12}y_{21}) \otimes_{F(Y)^H} Z' = 
UD_1'$ in the notation defined before 3.4. 

We now turn to the promised construction of a generic way 
to make a central simple algebra be in the image of restriction. 
Let $A/F$ be a central simple algebra of degree $r$ and suppose $K/F$ is a
separable degree $m$ 
field extension. Assume $B/K$ is also central simple of degree $r$. 
Form the generic splitting field $K(B \otimes_K A^{\circ})$ 
and define $F(A,B)$ to be the transfer 
$\Tr_{K/F}(K((B \otimes_K A^{\circ}))$ of this field to 
$F$. We set $KF(A,B)$ to be $K \otimes_F F(A,B)$. 

\proclaim Proposition 3.8. 
a) $B \otimes_K KF(A,B) \cong A \otimes_F KF(A,B)$. 
\smallskip
b) Suppose $B \cong B' \otimes_F K$ and $(r,m) = 1$. Then 
$F(A,B)(x_1, \ldots, x_{r^2-1})$ 
is purely transcendental over $F(A \otimes_F B'^{\circ})$. 
\smallskip
c) $F(A,B)/F$ has transcendence degree $m(r^2 - 1)$. 

\proof 
$K((A \otimes_F K) \otimes_K B^{\circ})$ is a subfield of 
$KF(A,B)$ and so $A \otimes_F KF(A,B)$ and $B \otimes_K KF(A,B)$ 
define equal elements in the Brauer group.  Having equal degrees, 
they are isomorphic. This proves a). As for b), 
form $L = F(A,B)(A \otimes_F B'^{\circ})$. Since $A$ and $B$ 
are equal in the Brauer group of $KF(A,B)$, taking corestrictions 
we have that $A^m$ and $B'^m$ are equal in the Brauer group of $F(A,B)$. 
Since $(r,m) = 1$, we have that $A$ and $B'$ are equal in the Brauer 
group of $F(A,B)$. Thus $L = F(A,B)(x_1, \ldots, x_{r^2 - 1})$. 
We can also write $L = F(A \otimes B'^{\circ})(A,B)$. 
For convenience, set $F' = F(A \otimes B'^{\circ})$. 
Then $A \otimes_F F'K \cong B \otimes_K F'K$, and so 
$F'(A,B)$ is the transfer of the rational field extension 
$F'K((A \otimes_F F'K) \otimes_{F'K} (B \otimes_{K} F'K)^{\circ})$ 
and thus is rational over $F'$. This proves b). 
The calculation of the transcendence degree is immediate.~\qed

The way to use the $F(A,B)$ construction to {\bf generically} 
force $B$ to be in the image of restriction is to make $A$ 
generic. That is, suppose $K/F$ has degree $m$ as above 
and $B/K$ is central simple of degree $r$. Let 
$Z = Z(F,r,2)$ be the center of the generic division algebra 
$UD = UD(F,r,2)$. Set $B_Z = B \otimes_K KZ$ and then set 
$F_R(B) = Z(UD, B_Z)$. That is, $F_R(B)$ 
is the field defined by generically forcing $B$ to be the image of 
$UD/Z$. Set $KF_R(B) = K \otimes_F F_R(B)$ 
and $UD_R = UD \otimes_Z F_R(B)$. We use 3.8 to show that: 

\proclaim Theorem 3.9. 
$B \otimes_K KF_R(B) \cong UD_R \otimes_{F_R(B)} KF_R(B)$. 
If $B = B' \otimes_F K$, then $F_R(B)(x_1,\ldots,x_{r^2-1})$ 
is rational over $Z(UD \otimes B'^{\circ})$ of degree $m(r^2 - 1)$
which is in turn rational over $F$ of transcendence degree 
$2r^2$. 
$F_R(B)$ has transcendence degree
$m(r^2 - 1) + r^2 + 1$ over $F$.

\proof The first two statements follow from 3.8 a) and b). 
By 3.1 $Z(UD \otimes B'^{\circ})$ is rational over $F$ 
of transcendence degree $2r^2$. The last statement follows 
by arithmetic. ~\qed

We apply this $F_R(B)$ construction some algebras that arose 
in section 2. For this reason, we again fix $r = 2$. Set 
$B = \Delta(F(Y_D)^{H'}/F(Y_D)^H,y_{12}y_{21})$ and  
$L' = F(Y_D)^H(y_1,\ldots,y_{2m})$. Set $L = 
F(Y_D)^W(y_1,\ldots,y_{2m})$. 
Part a) below will be proved by comparing the definition  
of $L_R(B)$ and 3.4.  

\proclaim Theorem 3.10. 
\smallskip
\item{a)} $L_R(B)\! =\! (F(Y_D)^W(y_1,\ldots,y_{2m}))_R(B)\! =\! F_m\! = 
\! F(P,Q,P_1,\ldots,P_m,Q_1,\ldots,Q_m)^{S_m}$.  
\smallskip
\item{b)} $L_R(B)$ is rational over $F(P,Q)$ and $F$. 
\smallskip
\item{c)} $(F(Y)^W)_R(B)(x_1,x_2,x_3)$ is rational over $Z(F(Y)^W,2,2)$. 

\proof 
We begin with a). Theorem 3.4 describes $F_m$ as a transfer, 
and $L_R(B)$ is defined as a transfer. The proof of 
a) then amounts to the verification that they are transfers 
of the same fields up to isomorphism. To prove b), 
all we need to remark is that 3.3 and the fact $F(P,Q)/F$ is rational 
finishes the argument. This later fact, that the center of generic 
$2 \times 2$ generic matrices is rational over $F$, is due to 
Procesi ([P1])]. 
Part c) is direct from 3.9.~\qed 
\bigskip
\leftline{Section four: Finishing the proof} 
\medskip
We are ready to prove the major theorem 4.2. 
We first dispose of the $n = 2$ case. 

\proclaim Lemma 4.1. $F(V)^{PO_2}$ and $F(V)^{PSp_2}$ 
are both rational over $F$. 

\proof 
Of course, 
$Z(F,2)$ is rational over $F$. If 
$A$ is Brauer equivalent to $UD(F,2) \otimes UD(F,2)$, 
then $A$ is split. Thus we are done by 1.1.~\qed

\proclaim Theorem 4.2. Suppose $n = 2m$, $m$ is odd. 
Let $V = M_n \oplus M_n$ be the representation 
of $PSp_n$ and $PO_n$ given by conjugation on each 
component. Then $F(V)^{PSp_n}$ and $F(V)^{PO_n}$ 
are rational over $F$. 

\proof 
As we said at the beginning, it suffices to consider 
$PSp_n$ and thus all the machinery of sections 
one through three. We can also assume $m \geq 3$. 
Form $(F(Y)^W)_R(B)$ and note that this can be written as 
$((F(Y_D))_R(B))_{\beta}(Y_0)^W$ where $\beta$ 
denotes the extension $0 \to Y_D \to Y \to Y_O \to 0$. 
Since $B$, suitably extended, is in the image of 
restriction, we have that $((F(Y_D))_R(B))_{\beta}(Y_0) = 
((F(Y_D))_R(B))(Y_0)$ with non-twisted $W$ action. 
The extension $0 \to Y_O \to Y_O' \to I_m \to 0$ is easily 
seen to be defined by the image, $\gamma_O$, of the canonical 
$\gamma \in H^2(W,Y)$. But $\gamma_O$ is trivial when restricted to 
$H$, so must have order dividing $m$. If we write 
$\gamma = \gamma_2 + \gamma_m$, where $\gamma_i$ has order $i$, 
then $\gamma_O$ must also be the image of $\gamma_m$. 
Let $B_m/F(Y)^W$ be the central simple algebra associated 
to $\gamma_m$ of degree $m$. 
Then $(F(Y)^W(B_m))_R(B) = ((F(Y)^W)_R(B))(B_m) = 
(((F(Y_D))_R(B))(Y_0)^W)(B_m) = ((F(Y_D))_R(B))(Y_0')^W$. 
Of course $2\gamma = 2\gamma_m$ and so  
we can set $B_m'$ to be the degree $m$ algebra associated 
with $2\gamma$ and derive that by Tregub's result (0.3 b))
$$\eqalign{
((F(Y_D))_R(B))(Y_0')^W = ((F(Y)^W)_R(B))(B_m) = \cr
((F(Y)^W)_R(B))(B_m') = ((F(Y)^W)(B_m'))_R(B).\enspace\cr}
\leqno{(2)}$$ 

We will analyze further both ends of (2). Beginning on the left, 
$Y_O'$ is a permutation $W$ lattice and so by 0.2 
$((F(Y_D))_R(B))(Y_0')^W$ is rational over $(F(Y_D)^W)_R(B)$ 
of degree the rank of $Y_O'$ which is $4m(m-1)$. 
Since $4m(m-1) \geq 2m$ for $m \geq 3$, this last field is rational 
over $F$ by 3.10 b). 

Turning to the right end of (2), $B \otimes_{F(Y_D)^H} F(Y)^H$ 
is in the image of restriction, namely it is the image of the 
quaternion algebra over $F(Y)^W$ associated to $\gamma_2$. 
Thus, by 3.9, $(((F(Y)^W)(B_m'))_R(B))(x_1,x_2,,x_3)$ is rational 
over $((F(Y)^W)(B_m')$ of degree $3m + 8$. We showed in 1.7 
that $F(V)^{PSp_n}$ is rational over $F(Y)^W(B_m')$ of degree 
$2m^2$. But $2m^2 > 3m + 8$ if $m \geq 3$. This proves 4.2.~\qed 

\bigskip
\leftline{Section Five: Four times an odd}
\medskip
In this section we will assume $n = 4m$, $n$ is odd, and show 

\proclaim Theorem 5.1.  
$F(V)^{PSp_n}$ and $F(V)^{PO_n}$ are stably rational. 

As before, it suffices to consider the former case. 

Now when it comes to stable rationality we can quickly simplify 
our situation. In [K], [Sc] it was shown that $Z(F,n)$ was stable 
isomorphic to $Z(F,4)Z(F,m)$. We gave a later proof 
of this in [S1]. Let $L$ be the field from [S1] 
rational over both $Z(F,n)$ and $Z(F,4)Z(F,m)$. 
In [S1] we also showed that $UD(F,n) \otimes_{Z(F,n)} L 
\cong (UD(F,4) \otimes_{Z(F,4)} L) \otimes_L (UD(F,m) \otimes_{Z(F,m)} L)$. 
That is, $UD(F,n)\! /Z(F,n)$ is stably isomorphic to 
$UD(F,4)$ $\otimes UD(F,m)$. 
 
Let $B$ be Brauer equivalent to $UD(F,n) \otimes UD(F,n)$. 
Write $B = B_2 \otimes B_m$ where $B_2$ has $2$ power degree 
and $B_m$ has odd degree. By e.g. [S1] p. 392, $L(B)$ 
is stably isomorphic to 
$L(B_2)(B_m^{\circ})$  
which is in turn obviously stably isomorphic to the join 
$(Z(F,4)(B_2)(Z(F,m)(B_m))$. Now $B_m$ is Brauer equivalent 
to $UD(F,m) \otimes UD(F,m)$ so $Z(F,m)(B_m)$ is stably 
isomorphic to $Z(F,m)(UD(F,m))$ (we do not need Tregub's 
stronger result here).  Of course, this last field 
is rational over $F$. Thus we have shown: 

\proclaim Lemma 5.2. To prove 5.1 
it suffices to show $F(V)^{PSp_4}$ is stably rational 
over $F$. 

{From} now on, then, we only deal with the $n = 4$ case. 
We know $F(V)^{PSp_4}$ is rational over $F(Y_2)^W$. 
In this case $W = A \rtimes \Z/2\Z$ has order 8. 
Let $\sigma$ generate the $\Z/2\Z$ part and $A = <\sigma_1> \oplus 
<\sigma_2>$ such that $\sigma(\sigma_i)\sigma = \sigma_{3-i}$. 
Let $C$ be the subgroup generated by $\sigma\sigma_1$ 
of order 4. In particular, $W$ can also be described as a dihedral 
group of order 8. Note that $H' = W \cap S_3$ is just 
$<\sigma_2>$ and $H = A \rtimes S_{m-1}$ is just $A$ in this case. 
Finally let $B$ be the subgroup of order 4 generated 
by $\sigma_1\sigma_2$ and $\sigma$. 

Set the lattice $Y' = Y_D' \oplus Y_O'$ just as above, and 
note that  
the restriction of $Y' \to I[W/H']$ to $Y_O'$ is onto. 
Also note that $Y_O' \cong \Z[W]$ (generated by $y_{31}$). 
We will be using 0.2 several times in the rest of this 
argument. 

If $Y_4$  is the kernel of $Y_O' \to I[W/H']$, 
then we have $0 \to Y_4 \to Y \to Y_D' \to 0$. 
By 0.2, it suffices to examine $F(Y_4)^W$. 
Note that $Y_4$ has rank 5. The canonical generator 
of $\Z[W]$ maps to $d_3 - d_1$. It follows that $1 + \sigma \in Y_4$. 
Thus $\Z[G/<\sigma>] \subset Y_4$. 
We need to find one more element in $Y_4$ and it is 
clear that it is $(1 - \sigma_1)(1 - \sigma_2)$. 

Of course what we really have to explore is the 
lattice $Y_2$.  
Recalling 1.3 we can take for $Y_2$ (here we are using 0.2 again) 
the preimage in $\Z[W]$ of $d_1 + d_2 - d_3 - d_4$. That is, 
$Y_2$ is generated by $Y_4$ and $(\sigma_1 + \sigma_2)$. 
All together, $Y_2$ is generated by $\Z[W/<\sigma>](1 + \sigma)$, 
$1 + \sigma_1\sigma_2$ and $\sigma_1 + \sigma_2$. 

Let $B$ be as above, and $M \subset Y_2$ the $B$ module 
generated by $\Z[B/<\sigma>](1+\sigma)$ and $1 + \sigma_1\sigma_2$. 
Then it is easy to see that $Y_2 \cong \Ind_B^W(M)$. 
It helps a bit to work with $M$. To simplify, 
set $B = <a> \oplus <b>$ where $b = \sigma$ and $a = \sigma_1\sigma_2$. 
Let $x \in M$ be $1 + b$ and $y = 1 + a$. Then $M$ is generated 
by $x,y$ subject to the relations that $b$ fixes $x$, $a$ 
fixes $y$ and $(1 + a)x = (1 + b)y$. 
Note that this last relation can be equivalently written as 
$(1 + ab)(x-y) = 0$. 

We define an embedding $M \to \Z[B] \oplus \Z$ by sending 
$x$ to $(1 + b, 1)$ and $y$ to $(b + ab,1)$. Computing 
the cokernel we have the exact sequence $0 \to M \to \Z[B] \oplus \Z \to 
\Z[B/<ab>] \to 0$. Inducing up to $W$ we have 
$0 \to Y_2 \to \Z[W] \oplus \Z[W/B] \to \Z[W/<\sigma_1\sigma_2\sigma>] 
\to 0$. We have almost shown: 

\proclaim Lemma 5.4. To prove 5.1 it suffices to show 
$F(U)^W/W$ is rational for some faithful $W$ representation 
$U$. 

\proof The above argument reduces us to considering 
$F(Z[W] \oplus \Z[W/B])^W$ which can also be written 
$F(U')^W$ for the corresponding permutation representation  
$U' = F[W] \oplus F[W/B]$. By Endo-Miyata's result (0.2), 
stably rationality for one faithful representation implies 
it for all.~\qed

Let $U$ be the dual of the two dimensional representation 
$Fv_1 + Fv_2$ where $\sigma_i(v_i) = -v_i$ 
$\sigma_i(v_j) = v_j$ if $i \not= j$, and $\sigma(v_i) = 
v_{3-i}$. Then we can view $F(U) = F(v_1,v_2)$ and it is 
clear that $F(U)^W = F(v_1^2 + v_2^2, v_1^2v_2^2)$ 
(Note also that $W$ is a reflection group on $U$ so we could 
quote Chevalley). Theorem 5.1. is proven. 
\bigskip 

\leftline{Appendix} 
\medskip
The initial version of this paper we assumed the ground field 
$F$ was algebraically closed of characteristic 0. In this 
appendix we gather the arguments that allow us to assume 
$F$ is any field. There are two key results in this section. 
The first, A.1, shows that the standard section argument applies 
in any characteristic, avoiding use of Zariski's Main Lemma. 
The second result, in A.3, avoids use of the proof of 
Bogomolov's No name Lemma. We also note that this whole 
section only deals with the group $PSp_n$ but it is clear 
that other, less elementary methods, could be employed 
for arbitrary reductive groups. 

If $V$ is a representation of $PSp_n$, 
we say $V$ is {\bf good} if there is an affine open 
subset $U_1 \subset V$ such that all points of 
$U_1$ have trivial stabilizer. Let $M^{-} \subset M_n(F)$ be the 
space of skew symmetric matrices and $\Delta \subset M^{-}$ 
the subspace of diagonal matrices. 
Set $V^{-} = M^{-} \oplus V$ and $V^{\Delta} = \Delta \oplus V$. 
If $x \in V^{-}$ or $y \in V^{\Delta}$ we set $\bar x \in M^{-}$ 
or $\bar y \in \Delta$ to be the projection on the 
first summand.  Let $T_{PSp} \subset N_{PSp} \subset PSp_n$ 
be the maximal torus and its normalizer described in section one. 
Set $W = N_{PSp}/T_{PSp}$ to be the Weyl group. 

\proclaim Theorem A.1. Let $V$ be a good representation 
and let $M^{-}$, $\Delta$, $V^{-}$, $V^{\Delta}$ and $N_{PSp}$ 
be as above.   
Then $F(M^{-} \oplus V)^{PSp_n} = F(\Delta \oplus V)^{N_{PSp}}$. 

\proof 
The action $PSp_n \times V^{-} \to V^{-}$ 
is the obvious one. We consider the variety $PSp_n \times (V^{\Delta})$. 
This variety has the 
$N_{PSp}$ action $n \cdot (g,x) = (gn^{-1},nx)$ and we will 
write the quotient as $PSp_n \times_N (V^{\Delta})$. 
Since $gx = (gn)(n^{-1}x)$, 
the $PSp_n$ action induces $\phi: PSp_n \times_N (V^{\Delta}) 
\to (V^{-})$. 

\proclaim Lemma A.2. $\phi$ restricts to an open 
immersion on some open subset of $PSp_n \times_N (V^{\Delta})$.  

\proof 
Let $X \subset PSp_n \times V^{-} \times V^{\Delta}$ 
be the closed reduced subvariety of $(g,x,y)$ such that 
$gy = x$. Clearly the projection $X \to PSp_n \times V^{\Delta}$ 
is an isomorphism. $X$ has an action by $N$ via 
$n \cdot (g,x,y) = (gn^{-1},x,ny)$ and a commuting action 
by $PSp_n$ given by $g' \cdot (g,x,y) = (g'g,g'x,y)$. 
The $N$ action on $X$ translates to the $N$ action 
on $PSp_n \times V^{\Delta}$. The ``action'' map 
$PSp_n \times V^{\Delta} \to V^{-}$ translates to the 
projection $\pi: X \to V^{-}$. Since the generic element 
of $M^{-}$ is diagonalizable, $\pi$ is dominant. 
Of course, $\pi$ induces $X/N \to V^{-}$. 

This last projection can be factored into 
$X \to V^{-} \times \Delta \to V^{-}$ where the 
first map, $\pi_1$, is $(g,x,y) \to (x,\bar y)$ and the 
second map, $\pi_2$, is the obvious projection. 
Note that $T_{PSp} \subset N_{PSp}$ acts trivially on $V^{-} \times \Delta$ 
so that $\pi_1$ factors through $X/T$. 
Let $Y$ be the closure of the image of 
$\pi_1$. Note that $W$ acts on $Y$ but as $W$ acts 
trivially on $V^{-}$ we have an induced $Y/W \to V^{-}$. 

The image of $\pi_1$ contains the dense subvariety 
defined by $p_{\bar x}(t) = \prod_i(t - \theta_i)$ 
and $d_{\bar x} \not= 0$ where $p_{\bar x}(t)$ is the characteristic 
polynomial of $\bar x$, $d_{\bar x}$ is the discriminant of 
$p_{\bar x}(t)$, and the $\theta_i$ are the (diagonal) entries 
of $\bar y$.  Thus the extension $F(V^{-}) \subset F(Y)$ amounts 
to adjoining roots of the generic $p_x(t)$ and so 
$F(Y)^W = F(V^{-})$. This implies there is an open 
subset $U_Y$ of $Y$ with $U_Y/W \to V^{-}$ an open 
immersion. 

Thus the lemma will be proven if we show 
$X/T_{PSp}  \to Y$ restricts to an open embedding. 
To show this, we will show that there is a rational section 
$s: Y \to X$. If $(x,\bar y) \in Y$, let $\theta_i$ 
be the entries of $\bar y$. Then $\theta_i1 - x$ 
is singular, and so $(\theta_i1 - x))(\theta_i1 - x)^* = 0$ 
where $z^*$ is the adjoint of $z$. Let $v_i$ be the first column 
of $(\theta_i1 - x)^*$. For some $U_Y' \subset Y$ open, all 
$\theta_i$ are distinct and all 
the $v_i \not= 0$. The $v_i$ are, of course, eigenvectors 
for $x$ with eigenvalue $\theta_i$. By the definition of 
the skew involution $J$, $\theta_{2i-1} + \theta_{2i} = 0$. 
Let $(v,v')$ be the skew form associated with $J$. 
That is, if $e_i$ is the standard basis, $(e_i,e_j) = 0$ 
if $|i-j| \not= 1$ and $(e_{2i-1},e_{2i}) = 1 = -(e_{21},e_{2i-1})$. 
Now it is immediate that $(v_i,v_j) = 0$ if $\theta_i + \theta_j 
\not= 0$. Since the form is nondegenerate, 
it is also immediate that $(v_{2i-1},v_{2i}) \not= 0$. 
We define $v_{2i-1}' = v_{2i-1}$ and 
$v_{2i}' = (1/(v_{2i-1},v_{2i}))v_{2i}$. If $a'(x,\bar y)$ 
is the matrix with $i$ column $v'_i$, then 
$a'(x,\bar y) \in Sp_n$  and $a'(x,\bar y){\bar x}a'(x,\bar y)^{-1} = 
\bar y$. Thus if $a(x,\bar y)$ is the image of $a'(x,\bar y)$ 
in $PSp_n$, the map 
$(x,\bar y) \to (a(x,\bar y), x,\bar y + a'(x,\bar y)(x - \bar x))$ 
is a section $U_Y' \to X$ of $X \to Y$. The composition 
$U_Y' \to X \to X/T_{PSp}$ is a rational inverse to the 
map $X/T \to Y$.  
All together, Lemma A.2 is proven.

We now turn to the proof of Theorem A.1. 
Let $U' \subset M^{-}$ be the affine open subset of elements $\bar x$ 
with $d_{\bar x} \not= 0$. Then  
$U'$ is a union of $PSp_n$ orbits, all of which are closed. 
Furthermore, any two elements of $U' \cap \Delta$ are in the 
same $N_{PSp}$ orbit. 
Let $p: V^{-1} \to M^{-}$ be the projection. 
There is an affine open $U'' \subset 
V^{-}$ such that all points of $U''$ have trivial stabilizer. 
For example, we can choose $U'' = M^{-} \times U_1$.  
  
Let $U \subset V^{-}$ be $p^{-1}(U') \cap U''$.  
Set $U_{\Delta} = V^{\Delta} \cap U$. Then restriction 
defines a morphism $F[U] \to F[U_{\Delta}]$ which obviously 
induces $\phi: F[U]^{PSp_n} \to F[U_{\Delta}]^{N_{PSp_n}}$. 
If $f \in F[U]^{PSp_n}$ satisfies $\phi(f) = 0$, 
then $f$ is 0 on $U_{\Delta}$. But all the 
$PSp_n$ orbits of $U$ meet $U_{\Delta}$, and $f$ is constant 
on such orbits, so $f = 0$. Thus $\phi$ is an injection. 
By 1.4, $F(V^{-})^{PSp_n}$ 
is the field of fractions of $F[V^{-}]^{PSp_n}$ 
and hence of $F[U]^{PSp_n}$. 
Thus $\phi$ induces an embedding  
$\phi': F(V^{-})^{PSp_n} \to F(V^{\Delta})^{N_{PSp}}$. 

As a $T_{PSp_n}$ module, $V$ is a direct sum 
of spaces of the form $Fv_{\tau}$ where $\eta(v_{\tau}) = 
\tau(\eta)v_{\tau}$ and $\tau$ is a character of $T$ (e.g. [K-T] p. 343). 
Since $V$ is good, it follows that the $\tau$ that appear 
generate the full group of characters of $T$. 
Since $N_{PSp}$ also acts on $T$, it follows that 
if $\tau$ appears in $V$ then so does $\tau^{-1}$. 
Thus 1.4 applies to $V^{\Delta}$ and $F(V^{\Delta})^{N_{PSp}}$ 
is the field of fractions of $F[V^{\Delta}]^{N_{PSp}}$. 

We can use Lemma A.2 to show $\phi$ is birational. 
It suffices to show that if $f \in F[V^{\Delta}]^{N_{PSp}}$, 
then $f$ is in the image of $\phi'$. 
Let $f' \in F[PSp_n \times V^{\Delta}]$ be defined by $f'(g,v) = 
f(v)$. Since $f$ is $N$ invariant, $f' \in F[PSp_n \times_N V^{\Delta}]$ 
and so by A.2 $f'$ defines a rational function $f''$ on $V^{-}$ and the 
definition of $f$ and $f'$ show that $f''$ is $PSp_n$ invariant. 
It is clear that $f$ is the image of $f''$. This proves 
A.1, the first of the two results we needed in this section.~\qed 

The second result we need is the rationality of 
$F(M_n(F) \oplus M_n(F))^{PSp_n}/F(M^{-} \oplus M_n(F))^{PSp_n}$. 
When $F$ is algebraically closed of characteristic 0, 
this is an immediate consequence of the argument in Bogomolov's 
``No-Name'' lemma. The point is that all $PSp_n$ representations 
are completely reducible. In order to handle general $F$, 
we use A.1. to reduce to the case of finite 
groups, and then we use the result of Endo-Miyata (0.2). To be precise, 
we will prove the following.  

\proclaim Theorem A.3. Suppose $V$ is a good representation 
as in A.1 and $M^{-}$, $V^{-}$ are also as described there. 
Let $V' = V \oplus V''$ where $V''$ is any $PSp_n$ representation. 
Write $V'^{-} = M^{-} \oplus V'$.  
Then $F(V'^{-})^{PSp_n}/F(V^{-})^{PSp_n}$ is rational. 

\proof 
First of all, it is clear that $V'$ is also a good representation. 
Thus the theorem reduces to showing 
$F(V'^{\Delta})^{N_{PSp}}/F(V^{\Delta})^{N_{PSp}}$ is rational, 
where $V^{\Delta}$ is as above and of course $V'^{\Delta} = 
\Delta \oplus V'$. We observed above that for each 
$\tau \in \Hom(T,F^*)$ there was a monomial 
$m_{\tau} \in F[V]$ such that $\eta(m_{\tau}) = \tau(\eta)^{-1}m_{\tau}$. 
$V''$ has a basis of $v_i$ such that $\eta(v_i) = 
\tau_i(\eta)v_i$ and we set $m_i = m_{\tau_i}$. 
Let $V_1$ be the $F$ span of the $m_iv_i$. 
It is clear that $F(V'^{\Delta}) = F(V^{\Delta} \oplus V_1)$, 
$T$ acts trivially on $V_1$, 
and that $F(V'^{\Delta})^T = F(V^{\Delta})^T(V_1)$. 
Now $V_1$ is not preserved by $W$ so we set 
$V_2 = F(V^{\Delta})^TV_1$. 

We claim $V_2$ 
is preserved by the action of $W$. 
For any $\tau \in \Hom(T, F^*)$, set $V_{\tau}'' = 
\{v \in V'' | \eta(v) = \tau(\eta)v $ all $\eta \in T \}$. 
Since for each $\tau$ we chose a unique $m_{\tau}$, 
$m_{\tau}V_{\tau}''$ is a subspace of $V_1$. 
It is clear that for $\sigma \in N_{PSp_n}$, 
$\sigma(V_{\tau}'') = V_{\sigma(\tau)}''$. 
Then 
$$\sigma(m_{\tau}V_{\tau}'') = 
(\sigma(m_{\tau})/m_{\sigma(\tau)})m_{\sigma(\tau)}V_{\sigma(\tau)}.$$
But $(\sigma(m_{\tau})/m_{\sigma(\tau)}) \in F(V^{\Delta})^T$ 
and so it is clear $N$ (i.e. $W$) preserves $V_2$. 

Of course $F(V^{\Delta})^T(V_1) = F(V^{\Delta})^T(V_2)$. 
$W$ acts faithfully on $F(V^{\Delta})^T$ and so 
by Endo-Miyata (0.2) $F(V^{\Delta})^T(V_2)^W$ 
is rational over $(F(V^{\Delta})^T)^W$. 
But $F(V^{\Delta})^T(V_2)^W = F(V'^{\Delta})^N$ and 
$(F(V^{\Delta})^T)^W = F(V^{\Delta})^N$.~\qed
\bigskip

\leftline{References}
\def\hangbox to #1 #2{\vskip1pt\hangindent #1\noindent \hbox to #1{#2}$\!\!$}
\def\refn#1{\hangbox to 30pt {#1\hfill}}

\medskip
\refn{[BS]}
Berele, A., and Saltman, D.J., 
{\it The centers of generic division algebras with involution}, 
Israel J. Math., (1) {\bf 63} (1988). 

\refn{[B]}
Brown, K., {\it Cohomology of Groups}, 
Springer-Verlag, New York/Heidelberg/Berlin 1982.

\refn{[EM]}
Endo, S. and Miyata, T., 
{\it Invariants of finite abelian groups}, 
J. Math. Soc. Japan, {\bf 25} (1973), 7--26.

\refn{[J]}
Jacobson, N., 
{\it Basic Algebra II}, 
W.H. Freeman \& Co., San Francisco, 1980.

\refn{[K]}
Katsylo, P.I., 
{\it Stable rationality of invariants of linear representations 
of the groups $PSL_6$ and $PSL_{12}$}, 
(Russian) Mat. Zametki {\bf 48} (1990), no. 2, 49--52; 
translation in Math. Notes {\bf 48} (1990), no. 1--2, 751--753 (1991). 

\refn{[K-T]}
Knus, M.A., Merkurjev, A., Rost, M., and Tignol, J.-P., 
{\it The Book of Involutions}, 
Amer. Math. Soc Providence RI, 1998.  

\refn{[LN]}
Saltman, D.J., 
{\it Lectures on Division Algebras}, 
Amer. Math. Soc., Providence, RI, 1999.

\refn{[P]}
Procesi, C., 
{\it Noncommutative affine rings}, 
Atti Acc. Naz. Lincei, S VII, v. VII, fo {\bf 6} (1967), 239--255. 

\refn{[P2]}
Procesi, C., 
{\it The invariant theory of $n \times n$ matrices}, 
Advances in Math. {\bf 19} (1976), no. 3, 306--381.  

\refn{[R]}
Rowen, L.H., 
{\it Polynomial identities in ring theory}, 
Academic Press, New York, London 1980.  

\refn{[S]}
Saltman, D.J., 
{\it Retract rational fields and cyclic Galois extensions}, 
Israel J. Math. {\bf 47} (1984), 165--215.

\refn{[S1]}
Saltman, D.J., 
{\it A note on generic division algebras}, 
Contemporary  Math. {\bf 130}, (1992).

\refn{[ST]}
Saltman, D.J., and Tignol, J. - P., 
{\it Generic algebras with involution of degree $8m$}, 
preprint. 

\refn{[Sc]}
Schofield, A., 
{\it Matrix invariants of composite size}, 
J. Algebra {\bf 147} (1992), no. 2, 345--349.

\refn{[T]}
Tregub, S.L., 
{\it Birational equivalence of Brauer-Severi manifolds}, 
Uspekhi Mat. Nauk. {\bf 46} (1991), 217--218; 
Russian Math. Surveys {\bf 46} (1991), no. 6, 229. 
 
\end